\magnification=1100

\font\fourteenbf=cmbx10 at 14pt

\def\R{\hbox{\bf R}}
\def\Z{\hbox{\bf Z}}
\def\N{\hbox{\bf N}}
\def\C{\hbox{\bf C}}
\def\P{\hbox{\bf P}}

\centerline{\fourteenbf The gluing construction for normally generic}
\centerline{\fourteenbf $J$-holomorphic curves}
\medskip
\centerline{Jean-Claude Sikorav}

\bigskip

\noindent{\bf Abstract.}\hskip3mm Let $C_0=[\Sigma_0,f_0]$ be a stable
$J$-holomorphic curve  in an almost complex manifold $(V,J)$. Under an assumption of normal
genericity, 
we show that $C_0$ has in the space of
homologous
stable $J$-curves of the same genus a locally Euclidean  neighbourhood of the expected dimension given by Riemann-Roch. 
In dimension $4$, the normal genericity 
condition is satisfied
if i) $\langle c_1(TV),C_i\rangle>0$ for each component of $C_0$ and ii) the only singularities of $C_0$ are  
ordinary 
double points. We also give an extension of these results for curves containing a given finite set. As an
application, we show that  every symplectic surface of degree $3$ in $\C\P^2$
is symplectically isotopic to an algebraic curve.

\bigskip
\noindent AMS 53C65, 30G20, 58F05

\vskip1cm
\noindent{\bf 1. Introduction and statement of the results}
\medskip
The gluing technique for $J$-holomorphic curves has been introduced by A.
Floer in his
definition of ``symplectic Floer homology'' [F]. It first appeared in the following form: given two (or more) $J$-holomorphic curves which intersect and satisfy a 
suitable
genericity condition, their union can be approximated by a continuous family of $J$-holomorphic 
curves, the convergence near the intersection points being like that of $\{xy=t\}\subset\C^2$ to 
$\{xy=0\}$ for $t\in\C^*$, $t\to 0$. The proof is similar to C. Taubes' gluing construction for 
instantons, ie one first constructs an approximately $J$-holomorphic curve and then one uses the implicit
function theorem to find a true $J$-holomorphic curve nearby. This is very well explained in the Appendix A
to the book [MS] of D. McDuff and D. Salamon.
\medskip
[MS] treated only the case where all curves are rational and the intersections take place between different
curves, in which case the deformed curve is also rational: cf. the model $\{xy=0\}\to\{xy=t\}$, now considered
in $\C\P^2$. Later it was realized that this gluing construction can be generalized to any $J$-curve parameterized by a 
Riemann surface with nodes (or {\it nodal}): in other words, the components 
may have arbitrary genus, and have self-intersection. The idea of parameterizing a $J$-curve by a nodal Riemann surface is due 
to M. Kontsevich [KM], who used it to reformulate
Gromov's compactness theorem in terms of {\it stable} $J$-holomorphic-maps, viewing it as an ``embedded'' version of Deligne-Mumford's compactification of
${\cal M}_g$.

\medskip
\noindent{\it Remark.}\hskip3mm We prefer to use the terminology ``stable $J$-holomorphic {\it curve}'', or simply ``stable
$J$-curve''.
\medskip
To formulate this generalized
gluing result, we introduce the space \ $\overline{\mkern -4mu{\cal M}}_g(V,J,A)$ of stable $J$-curves of genus $g$ in $V$, 
in the homology class $A\in H_2(V;\Z)$.  A point in this space is an
isomorphism class $[\Sigma_0,f_0]$, where
$\Sigma_0$ is a nodal Riemann surface of genus $g$ and $f_0:\Sigma_0\to V$ is a $J$-holomorphic map, with 
$(f_0)_*(\Sigma_0)=A$. The stability condition means that the group $\Gamma_0={\rm Aut}(\Sigma_0,f_0)$ is finite. If 
$\Gamma_0=\{\rm Id\}$, $f_0$ is called simple (or primitive).
We shall denote by ${\cal M}_g(V,J,A)$ the subspace obtained when $\Sigma_0$ is smooth.
\smallskip
We shall also need 
the extended spaces 
$$\overline{\hbox{\cal M}}(V)=\bigcup_{J\in{\cal J}(V)}\overline{\mkern-4mu{\cal M}}_g(V,J)\,\,,\,\,
\overline{\hbox{\cal M}}_g(V,A)
=\bigcup_{J\in{\cal J}(V)}\overline{\mkern-4mu{\cal M}}_g(V,J,A),$$
where ${\cal J}(V)$ is  the space of almost complex structures
on $V$. They are equipped with natural continuous projections \
$\pi:\overline{\hbox{\cal M}}_g(V)\,,\,\overline{\hbox{\cal M}}_g(V,A)\to{\cal J}(V).$
\medskip
\noindent{\it Remark.} Ideally, one would like to consider structures of class ${\cal C}^\infty$,
but for the Fredholm analysis and in particular the use of the Sard-Smale theorem, one has to work with
structures of class ${\cal C}^{N,\alpha}$ with $N$ very large (cf. [MS])).
\medskip

Let $C_0=[\Sigma_0,f_0]$  be an element of \ $\overline{\mkern -4mu{\cal M}}_g(V,J,A)$.
Linearizing the $J$-holomorphy equation $\overline\partial_Jf=0$ for maps $f:\Sigma_0\to V$ 
defined on the fixed Riemann surface $\Sigma_0$, 
one obtains the operator
$$D_{f_0}:\Gamma(f_0^*TV)\to\Omega^{0,1}(f_0^*TV),$$
which is of the form $\overline\partial+a$ (generalized Cauchy-Riemann). It is Fredholm (for suitable differetiability
classes), with complex index 
(half the real index) given by Riemann-Roch:
$${\rm ind}_{\C}(D_{f_0})=\langle c_1(TV),A\rangle+n(1-g).$$
Here $c_1(TV)\in H^2(V;\Z)$ is the first Chern class of $(TV,J)$. Note that the formal dimension $i(A,g)$ of 
\ $\overline{\mkern -4mu{\cal M}}_g(V,J,A)$ (or of ${\cal M}_g(V,J,A)$) over $\C$ is
$$i(A,g)=\langle c_1(TV),A\rangle+(n-3)(1-g)={\rm ind} (D_{f_0})+(3g-3),$$
the difference $3g-3$ being ${\rm dim}({\cal M}_g)-{\rm dim\ Aut}(\Sigma_0)$, accounting for the variation of the source 
as a 
Riemann surface and the equivalence relation.
\medskip
We shall say that $C_0=[\Sigma_0,f_0]$ is {\it parametrically generic} if $D_{f_0}$ is onto. 
We can now state the generalized gluing result. With different notations, it can be found in the papers
 [RT] of Y. Ruan and G. Tian (Theorem 6.1, special case $\nu_t=0$), [LT] of J. Li and G. Tian (Proposition 3.4, special case $S=0$) and [FO] of K. Fukaya and K. Ono ([FO], 
Theorem 12.9, special case $E_\sigma=0$). It is also implicit in the work [Sie1] of 
B. Siebert.
\bigskip
\noindent{\bf Theorem 0.}\hskip3mm {\it Let $C_0=[\Sigma_0,f_0]$ in \  $\overline{\mkern-4mu{\cal M}}_g(V,J,A)$ be 
parametrically generic. Then  there
is  a local homeomorphism 
$$\phi_{J}:(\overline{\mkern-4mu{\cal M}}_g(V,J,A),C_0)\to(\C^m\times(\C^{i(A,g)-m} /\Gamma_0)\,,\,(0,[0])),$$
where the finite group $\Gamma_0={\rm Aut}(\Sigma_0,f_0)$ acts linearly. Moreover, let us number the nodes $1,\cdots,m$, and let
$\tau_1,\cdots,\tau_m$ be the $m$  first complex coordinates of $\phi_J$. Then
the curves
$C$ such that $\tau_i(C)=0$ are exactly those 
for which the parameterizing surface $\Sigma$ keeps the
$i$-eth node of $\Sigma_0$.
\smallskip
Finally, there exists a local homeomorphism 
$$\widetilde\phi:(\overline{\hbox{\cal M}}_g(V,A),C_0)
\to(\C^m\times(\C^{i(A,g-m)} /\Gamma_0)\times{\cal J}(V)\,,\,(0,[0],J))$$ 
whose restrictions 
$\phi_{J'}$ all have the same properties as $\phi_J$.}
\bigskip
Actually, the papers [LT], [FO], [Sie1] also treat the case when $D_{f_0}$ is not onto (in particular for non-simple curves, where
this condition is no longer generic in $J$), by introducing a complement 
to the image and studying
an equation of the type $\overline\partial_Jf=g$ where $g$ belongs to this complement: the solutions have 
then
nothing to do with $J$-holomorphic curves, though of course counting them may be interesting from the point of view of
symplectic geometry (general definition of Gromov-Witten invariants or of Floer homology). Here we have a different point of view
since we are interested in $J$-holomorphic curves {\it per se}.
\bigskip
There is another
possible linearization of the equation $\overline\partial_Jf=0$, in which  one also allows the complex 
structure of $\Sigma_0$ to vary (cf. [IS]). Considering it as an almost complex structure $j$ on $\Sigma_0$, one differentiates
the equation $\overline\partial_Jf={1\over2}(df+J\,df\,j)$
 with respect to $j$ as well as to $f$. At the point $(j_0,f_0)$, we get an operator $$\eqalign{\widetilde D_{f_0}&:\Omega_{j_0}^{0,1}(T\Sigma_0)\times\Gamma(f_0^*TV)\to\Omega^{0,1}(f_0^*TV)\cr
&\widetilde D_{f_0}(v,\xi)={1\over2}J\,df_0(v)+D_{f_0}(\xi).\cr}$$

Let us call $C_0=[\Sigma_0,f_0]$ {\it normally generic} if $\widetilde D_{f_0}$ is 
onto. The terminology is due to the fact that, when it is possible to define a normal linearized 
$\overline\partial$-operator $D^{N_0}$ as in ([G], 2.1.C$_1$, [HLS], [IS]), the surjectivity of $\widetilde D_{f_0}$ is equivalent
to that of $D^N$: see section 4.

\medskip

\bigskip
The main result of this paper is that the gluing construction works for a normally generic curve.
\bigskip
\noindent{\bf Theorem 1.}\hskip3mm{\it $C_0=[\Sigma_0,f_0]$ be a stable $J$-curve which is normally generic. 
Then all the conclusions 
of Theorem 0 hold.}
\bigskip
In dimension $4$, the normal genericity can be ensured by homotopical assumptions ([G], 2.1.C$_1$, [HLS], [IS]).
We need here a slight generalization for a curve parametrized by a nodal surface.

\medskip
\noindent{\bf Proposition 1.}\hskip3mm {\it Let $(V,J)$ be an almost complex manifold of dimension $4$, 
and let \ $[\Sigma_0,f_0]=C_0$  be a stable $J$-curve  in \ $\overline{\mkern-4mu{\cal M}}_g(V,J)$. We assume that $f_0$ is an embedding with distinct
tangents near the nodes, and
 that the restriction $f_i$ to each component
$\Sigma_i$
satisfies
$$\langle f_i^*c_1(TV),\Sigma_i\rangle>|df_i^{-1}(0)|\,\,\,,\,\, i=1,\cdots,r.$$
Then
$C_0$ is normally generic.}
\medskip
The application of Theorem 1 to this case is most significant 
when $f_0$ is an embedding with distinct tangents at the nodes, since then one obtains a description of all geometrically 
close $J$-curves. A $J$-curve for 
which $f_0$ is an embedding with distinct tangents at the nodes will be called {\it nodal}. 
Combining Proposition 1 with Theorem 1 and the adjunction formula, we obtain

\bigskip
\noindent{\bf Corollary 1.} {\it Let $C_0$ in \ $\overline{\mkern-4mu{\cal M}}_g(V,J,A)$ be a nodal $J$-curve in dimension
\ $4$, with components  $C_1,\cdots,C_r$ in the image,
satisfying
$$\langle c_1(TV),C_i\rangle>0\,\,\,,\,\, i=1,\cdots,r.$$
Then there
are local homeomorphisms 
$$\eqalign{\phi_{J}:(\overline{\mkern-4mu{\cal M}}_g(V,J,A),C_0)&\to(\C^m\times\C^{d(A)-m}\,,\,(0,[0])),\cr
\widetilde\phi:(\overline{\hbox{\cal M}}_g(V,A),C_0)
&\to(\C^m\times\C^{d(A)-m} \times{\cal J}(V)\,,\,(0,[0],J))}$$
where
$$d(A)={1\over2}(A.A+\langle c_1(TV),A\rangle),$$ with properties as above.}

\bigskip
These results can be extended to the case when the curve is required to contain a given finite set $F\subset V$ in the image. 
We denote \ $\overline{\mkern-4mu{\cal M}}_g(V,J,A)$ the subset of such curves, and by $\widetilde D_{f_0,F}$ the restriction of 
$\widetilde D_{f_0}$
to the subspace
$$\{(v,\xi)\in\Omega^{0,1}(T\Sigma_0)\times\Gamma(f_0^*TV)\mid\xi=0\,\,{\rm on}\,\,f_0^{-1}(F)\}.$$
Then Theorem 1 generalizes to Theorem 1', and Proposition 1 
generalizes to Proposition 1' (Section 6). Corollary 1 then becomes

\medskip

\noindent{\bf Corollary 2.}\hskip3mm {\it Let $C_0$ in \ $\overline{\mkern-4mu{\cal M}}_g(V,J,A)$ be a nodal $J$-curve in dimension
$4$, with components  $C_1,\cdots,C_r$ in the image. Let $F$ be a finite subset of the smooth part of the image. Assume that
$$|F\cap C_i|<\langle c_1(TV),C_i\rangle\,\,\,,\,\,i=1,\cdots,r.$$
 Then  there are local homeomorphisms
$$\eqalign{\phi_J^F:(\overline{\mkern-3mu{\cal M}}_g(V,J,A;F),C_0)&\to(\C^m\times\C^{d(A)-|F\cap C_0|}\,,\,(0,[0]))\cr
\widetilde\phi_F:(\overline{\hbox{\cal M}}_g(V,A;F),C_0)&\to
(\C^m\times\C^{d(A)-|F|} \times{\cal J}(V)\,,\,(0,[0],J))\cr}$$ with properties as above.}

\bigskip
As an application of this last result, we give a partial result
on the problem of isotopy classes of symplectic surfaces in $\C\P^2$.
\medskip

\noindent{\bf Theorem 3.}\hskip3mm {\it A symplectic surface of degree $3$ in
$\C\P^2$ is symplectically  isotopic to an algebraic curve.}

\bigskip
\medskip
\noindent{\bf Organization of the paper.}\hskip3mm  In section 2, we recall the definition of stable $J$-holomorphic
curves and give a local description of these objects. 

In section 3  we prove Theorem 1, following the method of [MS]. 
For the continuity of the constructions with respect to $\vec t\in\C^m$, we follow [RT].

In section 4 we define (under some mild restrictions, and only in dimension $4$) the normal $\overline\partial$-operator
and prove that its surjectivity is equivalent to the normal genericity.

In section 5 we recall some important properties
in dimension $4$ and use them to prove Proposition 1 and thus Corollary 1. 

In section 6 we state and prove 
the results about curves containing a fixed finite subset: Theorem 1', Proposition 1' and Corollary 2. 

Finally,
 in section 7 we study the isotopy problem for symplectic surfaces and  and prove Theorem 3.

\vskip1cm
\noindent{\bf 2. The space of stable $J$-curves of genus $g$ in $V$}

\medskip
Let $(V,J)$ be an almost complex manifold.
\bigskip

\noindent{\bf 2.1. Definition.}\hskip3mm The space \ $\overline{\mkern-4mu{\cal M}}_g(V,J)$ of stable $J$-holomorphic curves of genus $g$ in $V$
is the set of isomorphism classes $[\Sigma,f]$, where
\medskip
\item{(i)} $\Sigma$ is a connected Riemann surface with nodes, of arithmetic genus $h^1(\Sigma,{\cal O}_\Sigma)=g$; recall that
this is also the
 genus of any connected smooth deformation, and is given by
$$g=\sum_{i=1}^r\, (\widetilde g_i-1) +m+1=\sum_{i=1}^r\, \widetilde g_i+m-r+1,$$ where $r$ is the number of components, 
$\widetilde g_1,\cdots,\widetilde g_r$ the genera of the normalizations
and $m$ the number of nodes
\smallskip
\item{(ii)}
$f:\Sigma\to V$ is a $J$-holomorphic map
\smallskip
\item{(iii)}
${\rm Aut}\,(\Sigma,f)$ is finite
(stability condition); equivalently, $f$ is
nonconstant on
each unstable component in the sense of Deligne-Mumford (ie rational
 and with at most $2$ points of the normalization lying
over the nodes, or of genus $1$ with no point over the nodes).
\medskip
\noindent If $A$ is an element of $H_2(V;\Z)$, then \ $\overline{\mkern-4mu{\cal M}}_g(V,J,A)$ is the subset of 
\ $\overline{\mkern-4mu{\cal M}}_g(V,J)$
defined by the condition $f_*([\Sigma])=A$.
\bigskip
We give now a local description of the topology of \ $\overline{\mkern-4mu{\cal M}}_g(V,J)$
near a point  $C_0=[\Sigma_0,f_0]$, following [FO] (sections 9 and 10) and [Sie1] (sections 2 and 3). 

\medskip
We define the topology by describing
a countable fundamental system of neighbourhoods of $[\Sigma_0,f_0]$.
For this, we choose an arbitrary Riemannian metric $\mu$ on $V$, so that the area
$${\rm area}_\mu(f)=\int_{\Sigma}|\Lambda^2df|$$ of any map $f:\Sigma\to V$ is defined. We also choose an arbitrary 
metric on $\Sigma_0$ and a sequence $\epsilon_i>0$, $\epsilon_i\to0$.
\medskip
Then by definition, the neighbourhood ${\cal N}_i$ ($i\in\N$) of $[\Sigma_0,f_0]$ is the set of all curves $[\Sigma,f]$ such
 that there exists a  continuous map $\kappa:\Sigma\to\Sigma_0$ with the following properties:
\smallskip
i) $\phi$ is a diffeomorphism over $\Sigma_0\!\setminus\! N_0$, where $N_0$
is the set of nodes, and  $\phi^{-1}(N_0)$ is a union of circles
\smallskip
ii)  $|\overline\partial\phi|\le\epsilon_i|\partial\phi_i|$ over 
$\Sigma_0\!\setminus\! {\cal U}_{\epsilon_i}$ 
where $ {\cal U}_{\epsilon_i}$ is a $\epsilon_i$-neighbourhood of $N_0$

\smallskip
iii) ${\rm dist}_{C^0}(f\circ\phi^{-1},f_0)\to0$ on  $\Sigma_0\setminus N_0$.

\medskip
\noindent{\bf Remarks} 
\medskip
1) Note that  if $V$ is a point then $\overline{\mkern-4mu{\cal M}}_g(V,J)=\overline{\mkern-4mu{\cal M}}_g$,
the Deligne-Mumford compactification.
\medskip
2) We have given a rather weak description of the convergence of $f_i$ to $f$. 
Actually, the proof of Gromov's compactness theorem shows that one will also have, for some $\epsilon'_i\to0$,
${\rm dist}_{C^1}(f\circ\phi^{-1},f_0)<\epsilon'_i$ on $\Sigma_0\setminus {\cal U}_{\epsilon'_i}$
and ${\rm area}_\mu(f_{|E_i})<\epsilon'_i$ for each component $E$ of ${\cal U}_{\epsilon'_i}$. In fact, we shall 
need the fact that, in a suitable sense, $f$ is close to $f_0$ in the $L_1^p$-topology for any $p\in]2,+\infty[$.
\medskip
3) One can generalize the definition of \ $\overline{\mkern-4mu{\cal M}}_g(V,J)$ to include marked points, thus obtaining
spaces \ $\overline{\mkern-4mu{\cal M}}_{g,m}(V,J)$. As in 1), one has \
$\overline{\mkern-4mu{\cal M}}_{g,m}(V,J)
=\overline{\mkern-4mu{\cal M}}_{g,m}$ \ if $V$ is a point.
\medskip
Clearly, the topology defined in this way does not depend on the choice of $(\epsilon_i)$ or $\mu$.
Less obvious, but true, is the fact that it is Hausdorff.
Another easy property is that the homology class $A=f_*([\Sigma])\in H_2(V;\Z)$ is locally constant on \
$\overline{\mkern-4mu{\cal M}}_g(V,J)$, so that one has a 
topological sum
$$\overline{\mkern-4mu{\cal M}}_g(V,J)=\coprod_{A\in H_2(V;{\hbox{\bf Z}})}\overline{\mkern-4mu{\cal M}}_g(V,J,A).$$

To state Gromov's compactness theorem, we recall the notion of tameness: an almost complex manifold $(V,J)$ equipped with a 
metric 
$\mu$ is {\it tame} if for every compact $K\subset V$ and every $C>0$, there is a compact $K'\subset V$ such that every 
$J$-curve
$C=f(\Sigma)$, $\Sigma$ connected, which meets $K$, is contained in $K'$. This is the case for instance if $V$ is compact, 
or $(V,J,\mu)$ has a compact
quotient, or more generally has a bounded geometry in a suitable sense.
\medskip
\noindent{\bf Theorem.}\hskip3mm {\it Let $(V,J,\mu)$ be a tame almost complex manifold. 
Then the area functional ${\rm area}_\mu$ defines a proper
map on \ \ $\overline{\mkern-4mu{\cal M}}_g(V,J)$.}
\bigskip
In practice, an important special case is the following.
\medskip
\noindent{\bf Corollary.}\hskip3mm {\it Let $(V,\omega,J)$ be a compact symplectic
manifold equipped with an almost complex structure which is $\omega$-positive. Then all the spaces \
$\overline{\mkern-4mu{\cal M}}_g(V,J,A)$, $A\in H_2(V;\Z)$, are compact.}

\bigskip
\noindent{\bf 2.2. Explicit description of $J$-curves close to $[\Sigma_0,f_0]$}
\medskip
Following [FO], [LT] and [Sie], we explicitly describe the curves $[\Sigma,f]$ which are
close to some fixed point $[\Sigma_0,f_0]$ in   \ $\overline{\mkern -4mu{\cal M}}_g(V,J)$.

\medskip
1) One defines, for $\vec t\in\C^n$ small enough, the differentiable surface 
$$\Sigma_{\vec t}=
\left((\Sigma_0\!\setminus\!{\cal U}_{\vec t})\cup{\cal A}_{\vec t}\right)_{\displaystyle/{\sim}},$$ where
$${\cal U}_{\vec t}=\bigcup_{i=1}^m\,\phi_i(\Delta_{|t_i|}\times\{0\}\cup\{0\}\times\Delta_{|t_i|}),$$
the 
$\phi_i$ being  holomorphic embeddings from $\Delta\times\{0\}\cup\{0\}\times\Delta$ to  $\Sigma_0$  whose images  are
disjoint 
neighbourhoods of the nodes, and
$${\cal A}_{\vec t}=\coprod_{i=1}^m\,\psi_{i}(A_{t_i})\,\,,\,\,A_t=\{(x,y)\in\Delta^2\mid xy=t\}.$$

\smallskip
\noindent
We identify
$\psi_{i}(x,y)$ with $\phi_{i}(x,0)$ for $|x|>|y|$, 
and with $\phi_i(0,y)$ for $|x|<|y|$.

\medskip

2) For every $v\in H^1(T\Sigma_0)$ small enough, one defines on  $\Sigma_{\vec t}$  the almost complex structure 
(necessarily integrable)
$j_{\vec t,v}$ as follows. First, for $\vec t=\vec 0$ one chooses a  lift $\sigma$ 
of the natural map $\Omega^{0,1}(T\Sigma_0)\to H^1(T\Sigma_0)$, which 
is equivariant under $\Gamma_0$ and such that each $\sigma(v)$ vanishes on ${\cal U}_1$. Then, denoting $j_0$ the almost complex structure on
$\Sigma_0$, one sets
 $$j_v=({\rm Id}-{j_0\sigma(v)\over2})\,j_0\,({\rm Id}-{j_0\sigma(v)\over2})^{-1}$$
so that the cohomology class of ${\partial j_v\over\partial v}_{|v=0}=\widetilde v$ is $v$ (Kodaira-Spencer map). 
Then there is a unique $j_{\vec t,v}$ such that $j_{\vec t,v}=j_v$ on $\Sigma_0\!\setminus\!{\cal U}_{\vec t}$ and 
$\psi_i^{*-1}(j_{\vec t,v})$ is standard on each
$A_{t_i}$. 
\medskip
3) One ``stabilizes'' $\Sigma_0$ by marking points on each unstable component, ie such that $2g_i+m_i<3$. We mark \  
$k_i=3-(2g_i+m_i)$ points on the smooth part of $\Sigma_i$. The total number of marked points is
$$k=\sum\, k_i=h^0(T\Sigma_0)=\hbox{\rm dim Aut}(\Sigma_0).$$
We obtain a stable 
marked curve $(\Sigma_0,\vec p_0)$, where $\vec p_0=(p_{0,1},\cdots,p_{0,k})$, which defines a point \ 
$[\Sigma_0,\vec p_0]\in \,\overline{\mkern-3mu{\cal M}}_{g,k}$.
By the stability hypothesis
on $(\Sigma_0,f_0)$, we can do this
 in such a way that $f_0$ is an embedding on a compact neighbourhood $N=\coprod N_j$ disjoint from the nodes. 
In particular, we 
can find disjoint
local transversal hypersurfaces $(H_j\subset V)$ to the images $f_0(N_j)$.
\medskip
 One also requires that no marked point lie in
${\cal U}_1$, and that ${\cal U}_1$ be invariant by the finite group
${\rm Aut}(\Sigma_0,\vec p_0)$, and that each $\sigma(v)$ vanishes  on $\vec p_0$.

\medskip
On $\Sigma_{\vec t}$ one marks $\vec p_{v,t}=\vec p_0$. This makes sense for $\vec t\in\Delta^m$
since these points are then in $\Sigma_0\!\setminus\!{\cal U}_{\vec t}$ which is canonically identified with 
$\Sigma_{\vec t}\!\setminus\!{\cal A}_{\vec t}$. 

\medskip
Then we have the
\medskip
\noindent{\bf Proposition.} (cf. [LT], Lemma 3.1) \hskip3mm{\it Let ${\cal N}$ be a small enough neighbourhood ${\cal N}$ of $[\Sigma_0,f_0]$
in   \ $\overline{\mkern -4mu{\cal M}}_g(V,J)$. Then there is a continuous map $$\psi:{\cal N}\to\C^m\times{ H^1(T\Sigma_0)\over{\rm Aut}(\Sigma_0,f_0)}$$
with the following property: if $\psi(C)=(\vec t,[v])$ then $C$ has a unique representative 
$((\Sigma_{\vec t},j_{\vec t,v});f))$ such that 
\smallskip
i) 
$f(p_{0,j})\in H_j$ $(\forall j)$
\smallskip

ii) $f$ is close to $f_0$ on 
$\Sigma_0\!\setminus\!{\cal U}_1$.}
\medskip
This is defined as follows:
\medskip
Given $C=[\Sigma,f]$ in \ $\overline{\mkern -4mu{\cal M}}_g(V,J)$ close enough to $[\Sigma_0,f_0]$, 
we mark $\vec p=(p_1,\cdots,p_k)$ on 
$\Sigma$ such that $f(p_j)\in H_j\,\,(\forall j)$; the choice of $\vec p$ is not unique, but this
will be discussed later. We obtain a point $[\Sigma,\vec p]$
in the Deligne-Mumford compactification $\overline{\mkern-3mu{\cal M}}_{g,k}$.  By the local structure of this space 
(cf. [HM]), it has a representative $((\Sigma_{\vec t},j_{\vec t,v});\vec p_0)$ where $(\vec t,v)$ is small; $\vec t$ is uniquely
defined, and $v$ is defined up to the action of ${\rm Aut}(\Sigma_0,\vec p_0)$ on $H^1(T\Sigma_0)$. Also, it is possible to choose $v$ such that
$C=[(\Sigma_{\vec t},j_{\vec t,v};f]$ with $f$ satisfying i) and ii) of the proposition. Note that 
$f$ is unique once $(\vec t,v)$ has been chosen. 

\medskip
There remains to  identify the different possible choices for $(\vec t,v)$, ie find all $(\vec t',v')$ close to $(\vec 0,0)$ such 
that there exists $f':\Sigma_{\vec t}\to V$, $(j_{\vec t',v'},J)$-holomorphic and satisfying

\smallskip
$[\Sigma_{\vec t',v'},f']=[\Sigma_{\vec t,v},f]$
\smallskip
$f'(p_{0,j})\in H_j$ $(\forall j)$

$f'_{|\Sigma_0\setminus{\cal U}_1}$ is close
to $(\vec 0,0,f_{0|\Sigma_0\setminus{\cal U}_1})$. 
\smallskip
\noindent We want to see that this is equivalent to $\vec t'=\vec t$, $v'=\gamma.v$
for some $\gamma\in{\rm Aut}(\Sigma_0,f_0)$. We shall prove only the direct implication, the converse being easy.

\medskip
 The above conditions are equivalent to the existence of a biholomorphism 
$\phi:(\Sigma_{\vec t},j_{\vec t,v})\to(\Sigma_{\vec t'},j_{\vec t',v'})$
such that $f'\circ\phi=f$, $\phi(p_{0,j})\in H_j$ ($\forall j$) and  $\phi$ is close on 
$\Sigma_0\!\setminus\!{\cal U}_1$ to an element 
$\gamma\in{\rm Aut}(\Sigma_0,f_0)$.

One easily proves 
that the subgroup $\widetilde\Gamma_0$ of ${\rm Aut}(\Sigma_0)$ generated by ${\rm Aut}(\Sigma_0,\vec p_0)$ and 
${\rm Aut}(\Sigma_0,f_0)$ is still finite. Thus we can impose that it leaves ${\cal U}_1$ invariant, thus each $\gamma\in\widetilde\Gamma_0$
will define an isomorphism $\gamma_{\vec t,v}$ from $(\Sigma_{\vec t},j_{\vec t,v})$
to $(\Sigma_{\vec t},j_{\vec t,\gamma.v})$.

\medskip

Thus $\psi=\gamma_{\vec t,v}\circ\phi^{-1}:\Sigma_{\vec t'}\to\Sigma_{\vec t}$ \ is a
$(j_{\vec t,v'},j_{\vec t',\gamma.v})$-biholomorphism
which is close to the identity on $\Sigma_0\!\setminus\!{\cal U}_1$, and $f'(\psi(p_{0,j}))\in H_j$ ($\forall j$). 
This implies 
$\psi(p_{0,j})=p_{0,j}$ ($\forall j$), thus $\vec t=\vec t'$ and $\psi={\rm Id}$. Thus also $v'=\gamma.v$, which finishes the proof.
\medskip
\noindent{\bf Remark.}\hskip3mm Condition iii) in the description of the topology implies that 
$f$ is close to $f_0$ on $\Sigma\setminus{\cal U}_{\vec t}$, that is on almost all its domain of definition.

\bigskip
\noindent{\bf 2.3. Decomposition according to the topology of the source}
\medskip
Denote by $R_g$ a set of representatives for each topological type of smooth surface with nodes 
of arithmetic genus $g$. For $\Sigma_0\in R_g$, we define ${\cal M}_{\Sigma_0}(V,J)$ as 
the subspace of \ $\overline{\mkern -4mu\cal M}_g(V,J)$ given by 
the curves $[\Sigma,f]$ with 
$\Sigma$
diffeomorphic to $\Sigma_0$. Then by definition we have a disjoint decomposition
$$\overline{\mkern -4mu{\cal M}}_g(V,J)=\bigcup_{\Sigma\in R_g}\,{\cal M}_{\Sigma}(V,J).$$
Let us call each element of this decomposition an {\it (equitopological) stratum}. 

\smallskip
One defines similarly
the strata ${\cal M}_{\Sigma,F}(V,J)$ where $(\Sigma,F)$ is a nodal surface and $F$ a finite subset 
of the smooth part. Thus if $R_{g,m}$ is a set of representatives of the topological types 
of $(\Sigma,F)$ for $\Sigma$ of genus $g$ and $|F|=m$, we have
 a disjoint decomposition
$$\overline{\mkern -4mu{\cal M}}_{g,m}(V,J)=\bigcup_{\Sigma\in R_{g,m}}\,{\cal M}_{\Sigma}(V,J).$$

\medskip
\noindent{\bf Local structure of \ ${\cal M}_{\Sigma_0}(V,J)$}
\medskip
1) If $\Sigma_0$ is smooth of genus $g$, ${\cal M}_{\Sigma_0}(V,J)$ is naturally homeomorphic to the space 
denoted \ ${\cal M}_g(V,J)$ by [IS]. Note that this last space has, under some genericity assumption, a structure
of smooth orbifold (cf. [IS] and section 4). The smooth structure comes from the fact that, 
up to a finite group, the space is naturally a subspace of a Banach manifold of maps.

\medskip
2) More generally, assume that $\Sigma_0$ is nodal, with components $\Sigma_1,\cdots,\Sigma_r$, 
but restrict to the case where ${\rm dim}(V)=4$ and  $f_0$ is an embedding with distinct tangents 
near each node. 
Then  ${\cal M}_{\Sigma}(V,J)$ is homeomorphic, and even diffeomorphic to an {\it open} subset in 
the product $\prod_{i=1}^r{\cal M}_{\widetilde\Sigma_i}(V,J)$, defined by suitable intersection properties.
\bigskip
\noindent{\it RemarK;}\hskip3mm The main problem in the study of spaces of $J$-curves is that, although equitopological strata are well 
understood, the way they ``hang  together'' in general is not.
\vskip1cm
\noindent{\bf 3. Proof of Theorem 1}
\bigskip
By 2.2, every $J$-curve $C$ close to $[\Sigma_0,f_0]$ is of the form $[\Sigma_{\vec t},j_{\vec t,v},f]$ with $\vec t\in\C^m$ well defined and
$v\in H^1(T\Sigma_0)$ determined
up to $\Gamma_0={\rm Aut}(\Sigma_0,f_0)$, $f$ being well defined once $v$ has been chosen. Also, $\vec t$ 
and $v$ tend to $0$ as $C\to C_0$, and $f\to f_0$ on the complement of the nodes (remark at the end of 2.2). Following [MS], [RT], 
[LT] and [FO], we shall find these curves by solving the equation for $(v,f)$ for a fixed small $\vec t$:
\medskip
1) First, we construct  a map $w_{\vec t}:\Sigma_{\vec t}\to V$ 
which is approximately $(j_{\vec t,0},J)$-holomorphic.
\smallskip
2) Then we write $f$ in the form $f=\exp_{w_{\vec t}}(\xi)$ where \ $\xi$ belongs to 
$\Gamma(\Sigma_{\vec t},w_{\vec t}^*TV\,{\rm rel}\,\vec p_0)$,
the space of sections of $w_{\vec t}^*TV$ which are of class $L_1^p$ for some $2<p<+\infty$ and are tangent to
$H_j$ at each point $p_{0,j}$.
The equation for $(v,f)$ 
can now be written in the form ${\cal F}_{\vec t}(v,\xi)=0$, where the target space of
${\cal F}_{\vec t}$ is $\Omega^{0,1}(w_{\vec t}^*TV)$, the space of forms of class $L_1^p$.
\smallskip
3) The equation is to be solved for $v$ small and $\xi$ small in the $L^\infty$-norm. We prove 
that in fact $\xi$ is necessarily 
small in any  norm $L_1^p$. We shall then work with such a norm on the source, with the corresponding 
$L^p$-norm on the target, so that ${\cal F}_{\vec t}$ becomes a Fredholm map.
\smallskip
4) We prove that ${\cal F}_{\vec t}$ satisfies the hypotheses of the implicit function theorem, 
with some uniformity in $\vec t$. The key point is that its differential
at zero is onto with a right inverse bounded independently of $\vec t$. For $\vec t=\vec 0$ it is exactly
the normal genericity. The existence of the inverse for $\vec t$ small is not obvious since the dependence of 
${\cal F}_{\vec t}$ is $\vec t$ is not quite continuous, but the construction of [MS] can be extended to our case.

\smallskip
5) A delicate point is the continuity in $\vec t$ of the solutions, for which we use an idea of [LT].
 
\bigskip
\noindent{\bf 3.1. Construction of an almost $J$-holomorphic map on $\Sigma_{\vec t}$}\hskip3mm \medskip
Let  $\vec t\in\C^m$ be such that $|\vec t]=\max |t_i|<1$. We construct a map
$w_{\vec t}:\Sigma_{\vec t}\to V$ (corresponding to $w_R$ in [MS])
as follows: $w_{\vec t}=f_0\circ\widetilde\rho_{\vec t}$, with
$$\left\{\eqalign{\widetilde\rho_{\vec t}&={\rm Id}\,\,{\rm on}\,\, \Sigma_{\vec t}\!\setminus\!{\cal A}_{\vec t}\cr
\widetilde\rho_{\vec t}(\psi_{i}(x,y))&=\phi_i(\rho(|t_i|^{-1/4}|x|)\,x,0)\,\,{\rm if}\,\,|x|\ge|y|\cr
\widetilde\rho_{\vec t}(\psi_{i}(x,y))&=\phi_i(0,\rho(|t_i|^{-1/4}|y|)\,y)\,\,{\rm if}\,\,|x|\le|y|.\cr}\right.$$ 
Here 
$\rho:\R^+\to[0,1]$ is a smooth cutoff function such that
$$\rho(s)=0\,\,{\rm for}\,\,s\le1\,\,{\rm and}\,\,\rho(s)=1\,\,{\rm for}\,\,s\ge2.$$
In particular:
$$w_{\vec t}(\psi_i(x,y))=f_0\circ\phi_i\,(0,0)\,\,{\rm if}\,\, |x|\,\,{\rm and}\,\,|y|\le|t_i|^{1/4}.$$

\medskip
\noindent{\it Comments.} We have replaced the $\delta$ of [MS], which 
was small but not infinitesimal, by $|t_i|^{1/4}$: this is to get the continuity in $\vec t$ 
of the constructions, notably that of the right inverse $R_{\vec t}$.
\bigskip
 To measure the almost-holomorphy of $w_{\vec t}$, we equip $\Sigma_{\vec t}$ with a metric which is fixed on
$\Sigma_0\!\setminus\!{\cal U}_1$
and induced on ${\cal A}_{\vec t}$ by the embeddings $A_{t_i}\subset\C^2$. 
\medskip

\noindent{\bf Proposition.}\hskip3mm{\it  One has for every $2<p<\infty$ the estimate
$$||\overline\partial_{j_{0,\vec t}}w_{\vec
t}||_{L^p}=O(|\vec t|^{1/2p}).$$}
\noindent{\bf Proof.} \hskip3mm We evaluate the contribution  to
$||\overline\partial_{j_{0,\vec t}}w_{\vec t}||_{L^p}^p$ of $A_{i,t_i}$, or
rather the
half $\{|x|\ge|t_i|^{1/2}\}$. Then
$\overline\partial_{j_{0,\vec t}} w_{\vec t}$ is $O(|t_i|^{-1/4}|x|)$ pointwise, and
vanishes except if
$|t_i|^{1/4}\le|x|\le2|t_i|^{1/4}$. Thus this contribution is at most of the
order of

$$(|t_i|^{-1/4})^p.\int_{|t_i|^{1/4}}^{2{|t_i|^{1/4}}}
\,r^{p+1}dr=O\big((|t_i|^{-p/4}).(|t_i|^{p/4+1/2})\big)=O(|t_i|^{1/2}).$$
Thus the total integral is $O(|\vec t|^{1/2})$, which proves (ii).

\bigskip
\noindent{\bf 3.2. Setting up the equation}
\medskip
We want to solve, for a given small $\vec t\in\C^m$, the equation
$$\overline\partial_{j_{\vec t,v},J}f:={1\over2}(df+J\,df\,j_{\vec t,v})=0$$ 
where the unknown is $(v,f)$ with \ $v\in H^1(T\Sigma_0)$ and $f\in{\rm Map}
(\Sigma_{\vec t},V)$ such that $f(p_{0,j})\in H_j\,\,(\forall j)$. Also, \ $v$ and
the ${\cal C}^0$ distance ${\rm dist}(f,w_{\vec t})$ should be small. 
\smallskip
Choosing the metric on $V$ so that the transversals $(H_j)$ are totally geodesic, we can write the unkown
$f$ uniquely
in the form $$f=\exp_{w_{\vec t}}(\xi)\,\,,\,\,\xi\in\Gamma(w_{\vec
t}^*TV{\rm rel}\,\vec p_0)\,\,\hbox{\rm small in the}\,\,L^\infty\,\,{\rm norm}.$$

\medskip
The expression $\overline\partial_{j_{\vec t,v}}f$ takes values in $\Omega^{0,1}(f^*TV)$. Let us fix a 
$J$-complex connection on $TV$ and use it to define the parallel transport along 
geodesics. We then get an isomorphism $\Pi_{\vec t}$ from $\Omega^{0,1}(f^*TV)$
to  $\Omega^{0,1}(w_{\vec t}^*TV)$. Thus we have to solve the equation 
${\cal
F}_{\vec t}(v,\xi)=0$, where $${\cal
F}_{\vec t}:H^1(T\Sigma_0)\times\Gamma(w_{\vec
t}^*TV{\rm rel}\,\vec p_0)\to\Omega^{0,1}(w_{\vec t}^*TV)$$ 
is defined by
$${\cal F}_{\vec
t}(v,\xi)=\Pi_{\vec t}(\overline\partial_{j_{\vec t,v},J}(\exp_{w_{\vec t}}(\xi))).$$ 
More precisely, we want to describe all 
the solutions such that $|v|$ and $||\xi||_{L^\infty}$ are small.
\medskip
The map ${\cal F}_{\vec t}$ \ is defined on a ball $B(\epsilon_0)=\{(v,\xi)\mid||(v,\xi)||<\epsilon_0\}$ 
independent of $\vec t$, and satisfies
$$||{\cal F}_{\vec t}(0,0)||=O(|\vec t|^{1/2p}).$$
Another property of ${\cal F}_{\vec t}$ is that it is of class ${\cal C}^\infty$. To see this, note that its pointwise value 
can be written
$$({\cal F}_{\vec t}(v,\xi))\,(z)=A(w_{\vec t}(z),v,\xi(z)).\nabla w_{\vec t}(z)+B(w_{\vec t}(z),v,\xi(z)).\nabla \xi(z),$$
where $A$ and $B$ are ${\cal C}^\infty$. This also implies that for every $k\in\N$, $||{\cal F}_{\vec t}||_{{\cal C}^k(B(\epsilon_0))}$  
is bounded, with a bound independent of $\vec t$.

\bigskip
\noindent{\bf 3.3. Smallness of $\xi$ in the $L_1^p$-norm}
\medskip
\noindent{\bf Proposition.}\hskip3mm {\it For every $\epsilon>0$ there exists $\eta>0$ with the following property.
Let $f=\exp_{w_{\vec t}}(\xi)$ be $(j_{v,\vec t},J)$-holomorphic, with
 $|v|$, $|\vec t|$ and $||\xi||_{L^\infty}$ less than $\eta$. Then $||\xi||_{L_1^p}$ is less than $\epsilon$.}
\medskip
\noindent{\bf Proof.}\hskip3mm 
Let $U$ be any neighbourhood 
of the nodes. By the elliptic apriori estimates, we can find $\eta$ such that 
$||\xi_{|\Sigma_0\setminus U}||_{L_1^p}$ is less than $\epsilon\over2$. Thus the proposition is reduced to a local 
property near the nodes. 
\smallskip
Let us work in a local chart $(V,v)\approx(\C^n,0)$ near the image of a node, $J$-holomorphic at $v$. 
Then the holomorphy equation for $f$ can be written $\overline\partial w+q(w).\partial w=0$ where $q$ 
is a map from $\C^n$ to $\overline{\rm End}_{\C}(\C^n)$. Also, replacing the coordinate $w$ by 
$\delta w$ with $\delta$ arbitrarily small, we can assume that $q$ is arbitrarily small in 
the ${\cal C}^1$-topology. 

\smallskip
Using a Riemannian metric on $V$ which is a multiple of the standard one in this chart,
the proposition is reduced  to the following lemma.
\medskip

\noindent{\bf Lemma.}\hskip3mm For every $\epsilon$ there exists $\eta$ with the following property.
Let $q$ be a map from the polydisk $\Delta^n\subset\C^n$ to $\overline{\rm End}_{\C}(\C^n)$ and let
$w$ and $w'$ be maps from $A_t$ to $\Delta^n$ which satisfy
the equations
$$\overline\partial w+q(w).\partial w=0=\overline\partial w'+q(w').\partial w'=0$$ and the inequalities 
 $$\max\,(|t|,||q||_{{\cal C}^1},||w||_{L_1^p},||w'-w||_{L^\infty})<\eta.$$
Then $||w'-w||_{L_1^p(A'_t)}<\epsilon$, where $A'_t$ is the intersection of $A_t$ with 
$\Delta_{1/2}\times\Delta_{1/2}$.
\medskip
\noindent{\bf Proof of the lemma}
\smallskip
We first construct a right inverse $P_t$ to $\overline\partial:L_1^p(A_t)\to L^p\Omega^{0,1}(A_t)$ 
which is uniformly bounded.
An element of $L^p\Omega^{0,1}(A_t)$ can be decomposed in two parts, one with support in the subannulus
$A_t^+=\{|x|\ge|y|\}$, the other in the subannulus $A_t^-=\{|y|\ge|x|\}$. It suffices to solve 
$\overline\partial f=\alpha$ for $\alpha$ with support in $A_t^+$. 
\smallskip
We write $\alpha=g(x)\overline{dx}$ where $g$ is a function of class $L^p$ defined  on 
$\Delta\setminus\Delta_{|t|^{1/2}}$, and we define $$P_t(\alpha)(x,y)=Pg(x)$$
using the standard inverse
$$Pg(x)=-{1\over\pi}\int_\Delta {g(z)\over z-x}\, d\sigma_z,\,\,x\in\C.$$
Recall that $\overline\partial(Pg)=g$ on $\Delta\setminus{\rm Supp}(g)$ and $0$ elsewhere. By Vekua 
(or Calderon-Zygmund),
 $P$ is bounded from $L^p(\Delta)$ to $L_1^p(\Delta)$. Since the projection $(x,y)\mapsto x$
from
$A_t$ to $\C$ is Lipschitz, this implies that $P_t$ is uniformly bounded.
\smallskip
Let $\xi=w'-w$, so that $||\xi||_{L^\infty}$ is arbitrarily small. Taking the difference between the equations and developing $q(w+\xi)-q(w)$ to the first order, we get
the linear equation
$$\overline\partial \xi+ q.\partial\xi+a.\xi=0$$
where we now consider $q$, as well as $a$, as a map defined on $A_t$. These maps  are 
arbitrarily small in the $L^\infty$-norm. 
\smallskip
We want to prove that then a solution of this equation satisfies 
$||q||_{L_1^p}\le C||q||_{L^\infty}$ with $C$ independent of $t$. This can be done in a way similar to [Sik1] p.171-172
which treats the case of a disk. One uses a smooth cutoff function $\rho:\Delta^2\to[0,1]$, 
which vanishes near the boundary and equals $1$ on $\Delta_{1/2}\times\Delta_{1/2}$. Then setting $\xi_1=\rho\xi$ one obtains
$$\overline\partial \xi_1+ q.\partial\xi_1=-(\overline\partial\rho+q.\partial\rho).\xi-a.\xi_1=g.$$
The right-hand side $g$ is $O(||\xi||)$ in the $L^\infty$-norm. It suffices to bound $||\xi||_{L_1^p}$
by $||\xi||_{L^\infty}$, ie to bound $||\overline\partial\xi||_{L^p}$ and $||\partial\xi||_{L^p}$.
\smallskip
Write $\xi_1=P_t(\overline\partial\xi_1)+h$ where $h$ is holomorphic, and set
 $T_t=\partial P_t$ so that $T_t$ is uniformly bounded in $L^p$ and
$\partial\xi_1=T_t(\overline\partial\xi_1)+\partial h.$
\smallskip
Note that the definition of $P_t$ implies that $||P_t(\overline\partial\xi_1)||_{L^2}\le||\xi_1||_{L^\infty}$, thus 
$||h||_{L^2}\le C||\xi||_{L^\infty}$. Since $h$ is holomorphic, one has $||h||_{L^\infty(A'_t)}\le C||h||_{L^2(A'_t)}$
and a fortiori
$||h||_{L^p(A'_t)}\le C||h||_{L^2(A'_t)}$
Thus
$$\eqalign{||\partial\xi_1||_{L^p(A_t)}=||\partial\xi_1||_{L^p(A'_t)}
&\le C(||\overline\partial\xi_1||_{L^p(A_t)}+||h||_{L^p(A'_t)})\cr
&\le C(||\overline\partial\xi_1||_{L^p}+C||\xi||_{L^\infty}).\cr}$$
The equation for $\xi$ implies then
$||\overline\partial\xi_1||_{L^p}(1-C\eta)\le C'||\xi||_{L^\infty}$. Thus for $\eta$ small enough we bound 
$||\overline\partial\xi_1||_{L^p}$ by $||\xi||_{L^\infty}$, and then also $||\partial\xi_1||_{L^p}$ 
via the estimate above. This concludes the proof of the lemma and thus of the proposition.

\bigskip
\noindent{\bf 3.4. Fredholm property and surjectivity of $d{\cal F}_{\vec 0}(0,0)$}
\medskip
\noindent{\bf Proposition.}\hskip3mm {\it 
\smallskip
i) The operator $D_{\vec t}=d{\cal F}_{\vec t}(0,0)$ is Fredholm of index $i(A,g)-m$. 
\smallskip
ii) The cokernel of $D_{\vec 0}$ is identified with that of $\widetilde D_{f_0}$. In particular,  $D_{\vec 0}$ is onto if and only if
 $[\Sigma_0,f_0]$ is normally generic.}
\medskip
\noindent{\bf Proof.} \hskip3mm First of all, let us observe that the linearization
operators $D_{f_0}$ and $\widetilde D_{f_0}$ can be defined when one replaces $f_0$ with any map $w$ of class $L_1^p$, giving
operators $D_w$ and $\widetilde D_w={1\over2}J\,dw\,\oplus 
D_w$ with the same properties.

\smallskip
i) The operator $D_{\vec t}$ is the restriction of $\widetilde D_{w_{\vec t}}$ to $\sigma(H^1(T\Sigma_0))\times\Gamma(w_{\vec
t}^*TV{\rm rel}\,\vec p_0)$. Since $D_{w_{\vec t}}$ is Fredholm and $H^1(T\Sigma_0)$ has finite dimension, 
this proves that $D_{\vec 0}$ 
is Fredholm. Its index is
$$\eqalign{{\rm ind}(D_{\vec 0})&=h^1(T\Sigma_0)+{\rm ind}(D_{f_0}{\rm rel}\,\vec p_0)\cr
&=h^1(T\Sigma_0)+{\rm ind}(D_{f_0})- h^0(T\Sigma_0)\cr
&={\rm ind}(D_{f_0})-\chi(T\Sigma_0).\cr}$$
On the other hand,
$$\eqalign{\chi(T\Sigma_0)&=\chi(T\widetilde\Sigma_0)-2m=\sum(3-3g_i)-2m\cr
&=3-3g+3m-2m=3-3g+m,\cr}$$
thus ${\rm ind}(D_{\vec t})={\rm ind}(D_{f_0})+3g-3-m$, qed.

\smallskip
ii) The operator $D_{\vec 0}$ is the restriction of $\widetilde D_{f_0}={1\over2}J\,df_0\,\oplus D_{f_0}$
to 
$\sigma(H^1(T\Sigma_0))\times\Gamma(f_0^*TV{\rm rel}\,\vec p_0)$. We want to prove that this does not diminish the image, ie
that

\smallskip
- restricting from $\Omega^{0,1}(T\Sigma_0)$ to $\sigma(H^1(T\Sigma_0))$: write
$\widetilde v=\sigma(v)+\overline\partial u$ with $u\in\Gamma(T\Sigma_0)$, thus 
$D_{\vec 0}(\widetilde v,0)=D_{\vec 0}(\sigma(v),{1\over2}ju)$
\smallskip
- restricting from $\Gamma(f_0^*TV)$ to $\Gamma(f_0^*TV{\rm rel}\,\vec p_0)$:  there exists
$X\in H^0(T\Sigma_0)$  such that $\xi(p_j)-df_0(X(p_j))\in TH_j$, thus $D_{\vec 0}(0,\xi)=D_{\vec 0}(0,\xi-df_0\,X)$.
\bigskip
\noindent{\bf 3.5. Existence of a uniformly bounded right inverse for $D_{\vec t}$}

\medskip
\noindent{\bf Proposition.} {\it For $|\vec t|$ small enough, $D_{\vec t}$ has a right inverse which is uniformly bounded. }
\medskip
\noindent{\bf Proof.}\hskip3mm One constructs a quasi-inverse $Q_{\vec t}(\eta)$ as in [MS], ie an operator which satisfies
$$||D_{\vec t}\circ Q_{\vec t}-{\rm Id}||_{L^p}=o(1)\,\,\,,\,\vec t\to\vec0.$$

For this, we shall need another almost $J$-holomorphic map $u_{\vec t}$, corresponding to
$(u_R,v_R)$ of [MS], this time defined on $\Sigma_0$: $u_{\vec t}=f_0\circ\widehat\rho_{\vec t}$, with
$$\left\{\eqalign{\widehat\rho_{\vec t}&={\rm Id}\,\,{\rm on}\,\,\Sigma_0\!\setminus\!{\cal U}_{\vec t}\cr
\widehat\rho_{\vec t}(\phi_{i}(x,0))&=\phi_i(\rho(|t_i|^{-1/4}|x|)\,x,0)\,\,{\rm if}\,\,|t|^{1/2}\le|x|\le1\cr
\widehat\rho_{\vec t}(\phi_{i}(0,y))&=\phi_i(0,\rho(|t_i|^{-1/4}|y|)\,y)\,\,{\rm if}\,\,|t|^{1/2}\le|y|\le1.\cr}\right.$$
To compare $u_{\vec t}$ with $w_{\vec t}$, define $\Gamma_{\vec t}$ to be the union of the circles (or points)
$$\Gamma_i=\psi_i(\{(x,y)\in\C^2\mid xy=t_i\,\,{\rm and}\,\,|x|=|y|=|t_i|^{1/2})\},$$
and $\Gamma_{\vec t}$ the union of the $\Gamma_i$. Let
$\pi_{\vec t}:\Sigma_{\vec t}\!\setminus\!\Gamma_{\vec
t}\to\Sigma_0$ be the map which is the identity on $\Sigma_{\vec t}\!\setminus\!{\cal A}_{\vec t}$, and
satisfies
$$\pi_{\vec t}\circ\psi_i(x,y)=\phi_i(x,0)\,\,{\rm or
}\,\,\phi_i(0,y)$$ depending whether $|x|>|y|$ or $|x|<|y|$. It is a
biholomorphism onto
its image $\Sigma_0\!\setminus\!{\cal D}_{\vec t}$.
\smallskip
Then we have $$w_{\vec t}=\left\{\eqalign{&u_{\vec t}\circ\pi_{\vec t}\,\,{\rm on}\,\,\Sigma_{\vec t}\!\setminus\!{\Gamma_{\vec t}}\cr
&\phi_i(0,0)\,\,{\rm on}\,\,\Gamma_i\cr}\right.$$

Replacing $w_{\vec t}$ by $u_{\vec t}$ in the definition of ${\cal F}_{\vec t}$, one  obtains a map
$${\cal G}_{\vec t}:\Omega^{0,1}(u_{\vec t}^*TV)\to\Gamma(u_{\vec t}^*TV).$$ 
We have ${\cal G}_{\vec 0}={\cal F}_{\vec 0}$, thus $d{\cal G}_{\vec0}(0,0)=d{\cal F}_{\vec 0}(0,0)$
is onto. It is easy to see that
$u_{\vec t}$ converges  to
$f_0$ in the
$L_1^p$ topology as
$\vec t\to\vec 0$. Thus $d{\cal G}_{\vec t}(0,0)$ has a right inverse
$\widetilde R_{\vec t}$ for $|\vec t|$ small enough, which is uniformly bounded and is continuous in $\vec t$.
\medskip
The operator $Q_{\vec t}$ is defined as the composition of the three following maps:
$$\matrix{&{\pi_{\vec t}}^{*-1}\cup 0&&\widetilde R_{\vec t}&&{\rm Id}\times
E_{\vec t}\cr
\Omega^{0,1}(w_{\vec
t}^*TV)&\longrightarrow&\Omega^{0,1}(u_{\vec
t}^*TV)&\longrightarrow&H\times\Gamma(u_{\vec
t}^*TV)&\longrightarrow&H\times\Gamma(w_{\vec t}^*TV)\cr}$$
The first map is  ${\pi_{\vec t}}^{*-1}$ followed by the extension
by zero from
$\Sigma_0\!\setminus\!{\cal D}_{\vec t}$ to $\Sigma_0$. The operator $E_{\vec t}$ is defined  by
$$(E_{\vec t}\xi)(z)=\left\{\matrix{&\xi(z)&{\rm if}\,\,
z\notin{\cal A}_{\vec t}\cr
&\!\!\!\xi\,(\phi_i(x,0))+\beta_{|t_i|}(x)
\,
\big(\xi\,(\phi_i(0,y))-\xi\,(z_i)\big)
&\!\!{\rm if}\,\,
z=\psi_i(x,y),|x|\ge|y|\cr
&\!\!\!\xi\,(\phi_i(0,y))
+\beta_{|t_i|}(y)\,\big(\xi\,(\phi_i(x,0))-
\xi\,(z_i)\big)
&\!\!{\rm if}\,\,
z=\psi_i(x,y),|y|\ge|x|.\cr}\right.$$
Here 
$\beta_{\delta}(z)=\rho({4\log |z|\over\log\delta}),$
where $\rho:\R\to[0,1]$ is the preceding cutoff function, satisfying $\rho(s)=0$ for $s\le1$, $\rho(s)=1$ for $s\ge2$.
One easily shows the following properties:
$$\displaylines{\beta_\delta(z)=1\,\,{\rm for}\,\,|z|\le\delta^{1/2}\,\,,\,\,  
\beta_\delta(z)=0\,\,{\rm for}\,\,|z|\ge\delta^{1/4}\cr
\int_{\C}|\nabla(\beta_\delta( z))|^p\,|z|^{p-2}d\sigma_z=o(1)\,\,{\rm
for}\,\,\delta\to0.\cr}$$
Note that the cutoff takes
place at points $\psi_i(x,y)$ with
$|t_i|^{1/2}\le|x|,|y|\le|t_i|^{1/4}$, where $w_{\vec t}$
is equal to
the constant $f_0\circ\phi_i(0,0)$.
\medskip
Clearly, $E_{\vec t}$ is  uniformly bounded for
$\vec t\to\vec 0$, and thus the same is true for  $Q_{\vec t}$ .
\medskip

\noindent{\bf Lemma.}\hskip3mm {\it One has the estimate}
$$||D_{w_{\vec
t}}(E_{\vec t}\xi)-D_{u_{\vec
t}}\xi||_{L^p}=o(||\xi||_{L_1^p})\,\,,\,\,\vec t\to\vec 0.$$

\noindent{\bf Proof.}
\hskip3mm
We need to estimate $$\alpha=D_{w_{\vec
t}}(E_{\vec t}\xi)-D_{u_{\vec
t}}\xi.$$  This form vanishes  except at points
$z=\psi_i(x,y)$ with $|x|,|y|\le|t_i|^{1/4}$.
We can assume $|x|\ge|y|$, thus $|t_i|^{1/2}\le|x|\le|t_i|^{1/4}$. 
Since  $w_{\vec t}(z)=u_{\vec t}(\phi_i(x,0))=v_i$,
$D_{w_{\vec t}}$ and $D_{u_{\vec t}}$ are both equal to the usual
$\overline\partial$-operator for maps into the complex vector space
$T_{v_i}V$, thus $$\alpha_z=\overline\partial(E_{\vec
t}\xi)\,(\psi_i(z))-\overline\partial\xi(\phi_i(x,0)).$$
Also,
$\overline\partial\xi(0,y)={\pi_{\vec t}^{*-1}}\eta(0,y)=0$ since
$|y|\le|t_i|^{1/2}$,
thus
$$\alpha_z
=-\overline\partial\big(\beta_{|t_i|}(x)\big)\,.\,
\big(\xi(\phi_i(0,t_ix^{-1}))
-\xi(z_i)\big).$$
Since the metric on $A_{t_i}$ is here equivalent to $|dx|^2$, we can work
in the
coordinate $x$. Using the  inequality $|t_ix^{-1}|\le|x|$ and the continuous injection
$L_1^p\to
{\cal C}^{1-2/p}$, we have the
pointwise estimate
$$|\alpha_z|
\le{\rm const.}||\xi||_{L_1^p}.|\overline\partial\big(\beta_{|t_i|^{1/4}}(x)\big)|.|x|^{1-2/p}$$
Thanks to the integral inequality on $\nabla\beta_\delta$, we
obtain the lemma.
\medskip
We can now prove the quasi-inverse property. Let
$\alpha=D_{\vec t}\circ Q_{\vec t}(\eta)-\eta$ for a given
$\eta\in\Omega^{0,1}(w_{\vec t}^*TV)$. 
Setting $(v,\xi)=\widetilde R_{\vec
t}\circ
\pi_{\vec t}^{*-1}(\eta)$, we have $\alpha_z=df_0.v+D_{w_{\vec
t}}(E_{\vec t}\xi)$ with $df_0.v+D_{u_{\vec
t}}\xi=0$. Thus $\alpha=D_{w_{\vec
t}}(E_{\vec t}\xi)-D_{u_{\vec
t}}\xi$, and the lemma gives $||\alpha||_{L^p}=o(||\xi||_{L_1^p})$ as $\vec t\to\vec 0$. Since
 $||\xi||_{L_1^p}=O(||\eta||_{L^p})$, this proves
$||D_{\vec t}\circ Q_{\vec t}-{\rm Id}||=o(1)$ as $\vec t\to\vec 0$.

\bigskip
\noindent{\bf 3.6. Isomorphism from $\ker D_{\vec t}$ to $\ker D_{\vec 0}$}
\medskip
\noindent{\bf Proposition.} {\it Let $\chi_{\vec t}:H^1(T\Sigma_0)\times\Gamma(w_{\vec t}^*TV)\to\ker D_{\vec 0}$ be defined by sending $(v,\xi)$ 
to $(v_0,\xi_0)$ which minimizes
$$|v_0-v|^2+\int_{\Sigma_0\setminus{\cal U}_1}|\xi_0-\xi|^2\,d\sigma.$$

Then the restriction of $\chi_{\vec t}$ to $\ker D_{\vec t}$
 is a linear isomorphism for $|\vec t|$ small enough, and moreover there is a uniform estimate
$$||\chi_{\vec t}(v,\xi)||_{L^2}\ge C^{-1}||(v,\xi)||_{L_1^p}.$$}
\medskip
\noindent{\bf Proof.} \hskip3mm
Arguing by contradiction, assume that there exists $\vec t_n\to\vec 0$ and $(v_n,\xi_n)\in\ker D_{\vec t_n}$ such that 
$||(v_n,\xi_n)||_{L_1^p}=1$ and $||\chi_{\vec t_n}(v_n,\xi_n)||_{L^2}=o(1)$.  Denote 
$\chi_{\vec t_n}(v_n,\xi_n)=(v'_n,\xi'_n)$. By definition, $(v_n-v'_n,\xi_n-\xi'_n)$ is orthogonal to $\ker D_{\vec t_n}$
on $\Sigma_0\!\setminus\!{\cal U}_1$.

\smallskip

There is a subsequence $(v_n,\xi_n)$ which converges away from the nodes to $(v,\xi)\in\ker D_{\vec 0}$. Since 
$(v'_n,\xi'_n)$ converges to $0$ in $L^2(\Sigma_0\!\setminus\!{\cal U}_1)$, $(v,\xi)$ is orthogonal to $\ker D_{\vec 0}$
on  $\Sigma_0\!\setminus\!{\cal U}_1$, thus $(v,\xi)=(0,0)$.  This implies that $|v_n|=o(1)$ and $\xi_n\to0$ away from the nodes,
uniformly in the ${\cal C}^\infty$ topology. Also $D_{w_{\vec t_n}}\xi_n=0$ on each annulus $A_{t_{n,i}}$.

\smallskip
On each annulus $A_t=A_{t_{n,i}}$ this equation takes the form
$$\overline\partial \xi+q(w).\partial\xi+(dq(w).\xi).\partial w=0$$
with $||q||_{L^\infty}$ and $||dq||_{L^\infty}$ are arbitrarily small. Also, $||w||_{L^\infty}$ 
and $||dw||_{L^\infty}$ are bounded, and
 is $\xi$ is arbitrarily ${\cal C}^\infty$
small away from $(0,0)$. In particular $||\xi_{|\partial A_{t_i}}||_{L^\infty}$ is arbitrarily small. 
\smallskip
We shall get a contradiction via the lemma in 3.3, provided we can prove a ``weak maximum principle'', namely the inequality

$$||\xi||_{L^\infty}\le C||\xi_{|\partial A_{t_i}}||_{L^\infty}.$$
To this effect, we set $\widetilde\xi=\xi+q(w).\xi$ (compare [Sik2], proof of Proposition 1). We have
$$\overline\partial\widetilde\xi=-(\overline\partial (q(w)).\xi-(dq(w).\xi).\partial w=A.\widetilde\xi,$$
where $A$ is complex-linear and is $o(1)$ in the $L^\infty$-norm.
\smallskip
As in [Sik2], the fact that $P_t$ is onto with a bounded right-inverse $L^p\to L_1^p$ gives a solution 
$\Phi:A_t\to{\rm G\ell}(\C^n)$ of the resolvant equation $\overline\partial\Phi=A.\Phi$. 
This solution lives in $L_1^p$ and satisfies $\max\,(||\Phi||_{L^\infty},||\Phi^{-1}||_{L^\infty})\le2$
for $||A||_{L^\infty}$ small enough.
Then writing $\widetilde\xi=\Phi.h$, we have that $h$ is holomorphic, thus $||h||$ satisfies the usual maximum
principle on $A_t$. Thus the weak maximum principle is satisfied with $C=4$, which finishes the proof of the proposition.

\bigskip
\noindent{\bf 3.7. End of the proof of Theorem 1}
\bigskip
We have proved that the map
$${\cal F}_{\vec t}:H^1(T\Sigma_0)\times\Gamma(w_{\vec t}^*TV{\rm rel}\,\vec p_0)\to\Omega^{0,1}(w_{\vec t}^*TV),$$
defined for
$|\vec t|<1$ on the ball $\{|(v,\xi)|<\epsilon_0\}$, has the following properties:
$$\left\{\eqalign{&{\cal F}_{\vec t}\,\,\hbox{\rm is of class}\,\, {\cal C}^2\,\,{\rm and}\,\,||{\cal F}_{\vec t}||_{{\cal C}^2}\le C_1\cr
&|{\cal F}_{\vec t}(0,0)|\le C_2|t|^{1/2p}\cr
&{\rm for}\,\,|\vec t|<\epsilon_1,\,\,D_{\vec 0}\,\,\hbox{\rm has a right inverse}\,\,R_{\vec t}\,\,
\hbox{\rm such that}\,\,|R_{\vec t}|\le C_3.\cr}\right.$$
Denote $$K_{\vec t}=\ker D_{\vec 0}\,\,,\,\,
K_{\vec t}(\epsilon)=\{(v,\xi)\in K_{\vec t}\mid|(v,\xi)|<\varepsilon\}.$$
By the implicit function theorem, we obtain 
a map
$$\Psi_{\vec t}=(\nu_{\vec t},\Phi_{\vec t}):K_{\vec t}(\epsilon)\to
H^1(T\Sigma_0)\times\Gamma(w_{\vec t}^*TV\,{\rm rel}\vec p_0)$$ 
defined for $|\vec t|$ and $\epsilon$ small enough,
with the
following property.
 If we set
 $$f_{\vec t,v,\xi}=\exp_{w_{\vec t}}(\Phi_{\vec t}(v,\xi)):\Sigma_{\vec t}\to V,$$ 
then $$\left\{\eqalign{&f_{\vec t,v,\xi}\,\,{is}\,\,(j_{\vec t,\nu_{\vec t}(\vec t,v)},J)\hbox{\rm-holomorphic
on}\,\, \Sigma_{\vec t}\cr
&f_{\vec t,v,\xi}(p_{0,j})\in H_j\,\,(\forall i)\cr
&{\rm dist}(f_{\vec t,v,\xi},w_{\vec t})<\varepsilon.\cr}\right.$$
Thus $[\Sigma_{\vec t,\nu_{\vec t}(v,\xi)},f_{\vec t,v,\xi}]$ 
is an element of the neigbourhood ${\cal N}_\epsilon$ of $[\Sigma_0,f_0]$

\medskip
Conversely, paragraph 3.3 and the uniqueness in the implicit function theorem
imply that, for some $\epsilon_1<\epsilon$ and for every $\vec t$ such that $|\vec t|<\epsilon_1$, the following property holds:
for every map $h$ defined on $\Sigma_{\vec t}$ 

$$\left\{\eqalign{&h\,\,{\rm is}\,\,(j_{v'},J)\hbox{\rm-holomorphic
on}\,\, \Sigma_{\vec t}\cr
&h(p_{0,j})\in H_j\,\,(\forall i)\cr
&{\rm
dist}(h,w_{\vec t})<\epsilon_1\cr}\right.$$
there exists $(v,\xi)\in K_{\vec t}(\epsilon)$ such that $h=f_{\vec t,v,\xi}$ and $v'=\nu_{\vec t}(v,\xi)$.
\bigskip
Let us now define $\phi_J:{\cal N}\to\C^m\times {\ker D_{\vec 0}\over{\rm Aut}(\Sigma_0,f_0)}$. 
We start with 
$\Phi(C)=({\vec t},[(v,f)])$. We write $f=\exp_{w_{\vec t}}\xi$, and we project
$(v,\xi)$ orthogonally on $(v_0,\xi_0)\in\ker D_{\vec 0}$, ie
$(v_0,\xi_0)=\chi_{\vec t}.(v,\xi).$
This is well defined up to the action of ${\rm Aut}(\Sigma_0,f_0)$. Then we set
$$\phi_J(C)=(\vec t,[v_0,\xi_0]).$$
This is clearly continuous, there remains to see that it is locally bijective.

\medskip
This finishes the proof of Theorem 1 for $J$ fixed. When $J$ varies, it suffices to apply the parametric version of 
the implicit function theorem.

\vskip1cm
\noindent {\bf 4. Normal genericity and the normal $\overline\partial$-operator in dimension $4$}
\bigskip
In this section we define the operator
$D^{N_0}$ presented in the Introduction. It has been defined by [IS] for all curves parameterized  by a smooth surface $\Sigma_0$, 
we shall need to extend the definition when $\Sigma_0$ has nodes. We shall do it only in dimension $4$, 
since in higher dimension there are some complications and anyway it is probably not very useful. Thus we 
assume in 
this section that $(V,J)$ is a $4$-dimensional almost-complex manifold.
\medskip
\noindent{\bf 4.1. The operator $D^{N_0}$}  
\medskip
Let $\Sigma_0$ be a nodal Riemann surface, and let $f_0:\Sigma_0\to V$ be a $J$-holomorphic map, with no constant component. 
\medskip
1)  If $\Sigma_0$ is smooth, following [IS] one defines the line bundle $N_0=
f_0^*TV/\overline{df_0(T\Sigma_0)}$, and $D^{N_0}:\Gamma(N_0)\to \Omega^{0,1}(N_0)$ is induced by 
$D_{f_0}$. It is proved in [IS] (paragraphs 1.3 and 1.5, see also 5[B], paragraph 1.1) that $D^{N_0}$ is 
an operator of type $\overline\partial+a$ defined from $L_1^p$ sections to $L^p$ forms, and that
one has an exact diagram
$$\matrix{0\to&\Gamma(T\Sigma_0\otimes L)
&\longrightarrow&\Gamma(f_0^*TV)&\longrightarrow&\Gamma(N_0)&\to0\cr
&\overline\partial\downarrow&&D_{f_0}\downarrow&&D^{N_0}\downarrow\cr
0\to&\Omega^{0,1}(T\Sigma_0\otimes L)&
\longrightarrow&\Omega^{0,1}(f_0^*TV)&\longrightarrow&\Omega^{0,1}(N_0)&\to0\cr}\leqno(4.1)$$
Here $L=L(df_0^{-1}(0))$ is the divisor of zeros of \ $df_0$, counted with multiplicities. 
In particular, if $f_0$ is an immersion:
$$\matrix{0\to&\Gamma(T\Sigma_0)
&\longrightarrow&\Gamma(f_0^*TV)&\longrightarrow&\Gamma(N_0)&\to0\cr
&\overline\partial\downarrow&&D_{f_0}\downarrow&&D^{N_0}\downarrow\cr
0\to&\Omega^{0,1}(T\Sigma_0)&
\longrightarrow&\Omega^{0,1}(f_0^*TV)&\longrightarrow&\Omega^{0,1}(N_0)&\to0\cr}.$$

\medskip
2) If $\Sigma_0$ has nodes, we make the following assumption:
\medskip 
\item{(*)} {\it $f_0$ is an embedding near each node, with distinct tangents.}
\medskip
Consider the normalization
$\nu:\widetilde\Sigma_0\to\Sigma_0$ and the induced map $\widetilde f_0=f_0\circ\nu$. One can
associate
as in 1) the normal bundle $\widetilde N_0$ over 
$\widetilde\Sigma_0=\coprod_{i=1}^r\widetilde\Sigma_i$ and the operator $D^{\widetilde N_0}$. By definition, 
we set $N_0=\widetilde N_0$ and $D^{N_0}=D^{\widetilde N_0}$.
\medskip
Let us prove that the diagram (4.1) ``over $\widetilde\Sigma_0$'' remains exact ``over $\Sigma_0$''.
Only the first line is changed: 
$\Gamma(T\Sigma_0\otimes L)$ (resp. $\Gamma(f_0^*TV)$) can be identified 
with elements of  \ $\Gamma(\widetilde T\Sigma_0\otimes L)$ (resp. $\Gamma(\widetilde f_0^*TV)$)
vanishing on the inverse images $z^+$ and $z^-$ of any node (resp. having the same value on 
 $z^+$ and $z^-$).
\smallskip

The problem is to prove that $\Gamma(f_0^*TV)\to\Gamma(N_0)$ is still onto. This follows immediately 
from the fact that,
for each double point $v=\widetilde f_0(z^+)=\widetilde f_0(z^-)$, the natural map
$$T_vV\to{T_vV\over d\widetilde f_0(T_{z^+})}\oplus{TvV\over d\widetilde f_0(T_{z^-})}$$
is onto by the assumption (*).
\medskip
\noindent{\it Remark.} This is where we use the dimension $4$: in higher dimension, we would have to restrict
$D^{\widetilde N_0}$ to some suitable subspace of $\Gamma(\widetilde N_0)$.

\bigskip
\noindent{\bf 4.2. Index computation}
\medskip
Using  the exact diagram (4.1) and the index formula for $D_{f_0}$ p.2, one gets the following result. 

\medskip
\noindent{\bf Proposition.}
\hskip3mm{\it One has $${\rm ind}(D^{N_0})=\langle c_1(TV),A\rangle+g-1-m-|df_0^{-1}(0)|.$$}
\noindent Other proof: applying  Riemann-Roch separately to each component, one has 
${\rm ind}(D^{N_0})=$ \break $\sum_i\,(c_1(\widetilde N_i)+1-g_i)$. The short  exact sequence 
$T\widetilde\Sigma_i\to \widetilde f_0^*TV\to \widetilde N_i$ gives
$$c_1(\widetilde N_i)=\langle c_1(TV),A_i\rangle-2(1-g_i)-|df_i^{-1}(0)|.$$
Using the fact that $\sum\,(1-g_i)=1-g+m$, we get the proposition.
\medskip
\noindent{\it Remark.}\hskip3mm Note that $\langle c_1(TV),A\rangle+g-1=i(A,g)$.
\bigskip

\noindent{\bf 4.3. Normal genericity and surjectivity of $D^{N_0}$}
\medskip
Let $f_0:\Sigma_0\to V$ be a $J$-holomorphic map such that $D^{N_0}$ can be defined. Then 
the exact diagram (4.1) implies the isomorphism

$${\Omega^{0,1}(f_0^*TV)\over df_0(\Omega^{0,1}(T\Sigma_0\otimes L)\oplus D_{f_0}(\Gamma(f_0^*TV))}
\simeq\,{\rm coker}\,D^{N_0}.$$

Since $\widetilde D_{f_0}={1\over2}J\,df_0\oplus D_{f_0}$ and $J\,df_0=df_0\,j_0$, one has 
$${\rm coker}\,\widetilde D_{f_0}
={\Omega^{0,1}(f_0^*TV)\over df_0(\Omega^{0,1}(T\Sigma_0))\oplus D_{f_0}(\Gamma(f_0^*TV))}.$$
Since $\Omega^{0,1}(T\Sigma_0)\subset\Omega^{0,1}(T\Sigma_0\otimes L)$, there is a natural embedding
from ${\rm coker}\,D^{N_0}$ to ${\rm coker}\,\widetilde D_{f_0}$.
Furthermore, the induced map 
$H^{0,1}(T\Sigma_0)\to\Omega^{0,1}(T\Sigma_0\otimes L)$ is onto, thus if 
$\alpha\in\Omega^{0,1}(T\Sigma_0\otimes L)$
there exist $\beta\in\Omega^{0,1}(T\Sigma_0)$ and $\xi\in\Gamma(T\Sigma_0\otimes L)$ with 
$\alpha=\beta+\overline\partial\xi$, so that $df_0(\alpha)=df_0(\beta)+D_{f_0}\,df_0(\xi)$ belongs to ${\rm im}\,\widetilde D_{f_0}$.
\smallskip
Thus ${\rm coker}\,D^{N_0}={\rm coker}\,\widetilde D_{f_0}$, thus normal genericity is equivalent to the surjectivity of $D^{N_0}$.

\vskip10mm
\noindent{\bf 5. The case of dimension $4$}
\bigskip
\noindent{\bf 5.1. Adjunction formula}
\medskip
The most important special property for $J$-curves in dimension $4$ is the  
``positivity of intersections'' and in particular its corollary the ``adjunction formula'' (cf. [McD], [MW], [Sik2]). We state here a version allowing the source of the 
curves to have nodes.
\medskip
\noindent{\bf Proposition.} \hskip3mm{\it Let $f:\Sigma\to (V^4,J)$ be a simple
$J$-holomorphic map, where $\Sigma$ is a connected Riemann surface with $m$ nodes. Then to each
singularity $s$ of the image $S=f(\Sigma)$ one can associate a strictly
positive integer $\delta(s)$, with the property that
$$\sum_{s\in{\rm Sing}(S)}\,\delta(s)-m=g_a(A)-g_a(\Sigma).$$
Here $g_a(A)$ is the arithmetic genus of $A$, which is equal to ${1\over2}(A.A-\langle c_1(TV),A\rangle)+1$.
}
\medskip
Note that each node $z$ contributes at least $1$ to the sum on the left, and exactly $1$ if and only if it is embedded
with distinct tangents and $f^{-1}(f(z))=\{z\}$.
\bigskip
\noindent{\bf 5.2. Nodal curves}
\medskip
Let $(V,J)$ be an almost complex manifold of dimension $4$. We say that a $J$-holomorphic curve 
$C_0=[\Sigma_0,f_0]$
is {\it nodal} if $f_0$ is an embedding. By 5.1, the arithmetic genera $g_a(\Sigma)$ and $g_a(A)$ coincide, where
$A=f_*(\Sigma)$.
\smallskip

The image $f_0(\Sigma_0)=S$ is then a closed immersed real
surface
whose singularities are ordinary double points (or nodes), and whose tangent
bundle is
$J$-invariant: $JTC_0=TC_0.$ 
\smallskip
Conversely, 
using the integrability of almost complex structures
on surfaces, such a surface is of the form $S=f_0(\Sigma_0)$, where $[\Sigma_0,f_0]$ is nodal. 
Moreover, the pair
$(\Sigma_0,f_0)$ is determined by $S$ up to isomorphism, and thus $S$ determines a stable
curve $[\Sigma_0,f_0]$ in $\overline{\mkern -4mu{\cal M}}_g(V,J,A)$
where $A=[S]$ and $g=g_a(A)$.
\medskip
Note that 
$$i(A,g_a(A))=\langle c_1(TV,A)\rangle+g-1={1\over2}(A.A+\langle c_1(TV),A\rangle).$$
We shall denote this number by $d(A)$.
Note that
$${\rm ind}\,D^{N_0}=\sum_i\,{\rm ind}\,D^{N_i}=\sum_i\,i(A_i,g_i)=\langle c_1(TV),A\rangle+g-1-m.$$ If $f_0$ is an embedding, this is also
${1\over2}(A.A+\langle c_1(TV),A\rangle)-m$.
\bigskip
\noindent{\bf Proposition.}\hskip3mm {\it If \ ${\rm dim}(V)=4$ and $[\Sigma_0,f_0]$ is a 
$J$-holomorphic embedding, then 
every simple $J$-curve $[\Sigma,f]$ in \ $\overline{\mkern -4mu{\cal M}}_g(V,J,A)$ is also an embedding. In particular, this applies
to every curve
 sufficiently close to $[\Sigma_0,f_0]$. }

\medskip
\noindent{\it Remark.} Presumably, the same result holds in all dimensions,
but I do
not know how to prove it in the nonintegrable case. 
\medskip
\noindent{\bf Proof.} The adjunction formula above implies that $$\sum_{s\in{\rm Sing}(f(\Sigma))}=
\#(\hbox{\rm nodes of }\,\Sigma).$$
 Thus necessarily
each node is
embedded, and there is no other singularity on $C$ which means that $f$ is
an embedding.

\bigskip
\noindent{\bf 5.3. Automatic regularity}
\medskip
Another aspect of the positivity of intersections is the ``automatic regularity'' of spaces of $J$-curves
under homotopic conditions [G, 2.1.C$_1$] (cf. also [HLS]). We begin by recalling the linear version.

\medskip
\noindent{\bf Proposition.}\hskip3mm {\it Let $L$ be a complex line bundle over a smooth Rieman, surface of genus $g$, equipped with an operator
$D:\Gamma(L)\to\Omega^{0,1}(L)$ of the type $\overline\partial+a$. Then $D$ is onto provide $c_1(L)>2g-2$. }
\medskip
Now let $f_0:\Sigma_0\to (V^4,J)$ be a $J$-holomorphic map defined on a smooth surface of genus $g$.
When applying the proposition to $L=N_0$, the condition $c_1(L)>2g-2$ becomes $\langle c_1(TV),A\rangle>0$ 
if $f_0$
is an immersion, and $c_1(f_0^*TV),\Sigma_i\rangle>|df_0^{-1}(0)|$ if $f_0$ is nonconstant [IS].
\smallskip
More generally, assume that
 $\Sigma_0$ is nodal, with components $\Sigma_1,\cdots,\Sigma_r$, and  
$f_0$ is an embedding with distinct tangents 
near each node. Then $D^{N_0}$ is isomorphic to the product
of the $D^{\widetilde N_i}$. Using the formula for $C_1(\widetilde N_i)$ in 4.2, we get the

\medskip
\noindent{\bf Proposition 1.} {\it Assume that ${\rm dim}(V)=4$, $\Sigma_0=\cup_i\,\Sigma_i$ is nodal of genus $g$ with 
$m$ nodes, 
that $f_0$ is an embedding near the nodes and has no constant component, and 
$$\langle c_1(f_i^*TV),\Sigma_i\rangle>|df_0^{-1}(0)|\hskip3mm\forall i=1,\cdots,r.$$ Then $D^{N_0}$ is onto.}

\bigskip

\noindent{\bf Remark.}\hskip3mm The existence of a $J$-holomorphic curve such that $\langle c_1(TV),C_0\rangle>0$ and $C_0$ is not exceptional
on  a closed almost complex
 $4$-dimensional manifold
tamed by a symplectic structure $\omega$ is quite restrictive: by a theorem of Li and Liu [LL], this 
implies that $(V,\omega)$ is rational or ruled.

\bigskip
\bigskip
\noindent{\bf 6. Fixing points}
\medskip
The results that we have obtained can be generalized to curves containing a given finite subset 
$F$ in the smooth part. At the linearized level, it means studying the restriction of $D^{N_0}$ to
 the subspace of sections of $N_0$ which vanish on $\widetilde F$. 

\smallskip
The way to do this is explained
in [B], Lemma 4. One replaces $N_0$ by $N_0\otimes L(F)$, where $L(F)$ is the bundle associated to the divisor $F$
on $\Sigma_0$. Then $D^{N_0}$ induces an operator 
$$D^{N_0}_F:L_1^p\Gamma(N_0\otimes L(F))\to L^p\Omega^{0,1}(N_0\otimes L(F)).$$
The regularity $(L_1^p,L^p)$ is essential since $D_F^{N_0}$ contains terms of the type ${\overline z\over z}$
in a local coordinate near any point of $F$. Actually, one does this for the restriction
$N_i$ of $N_0$ to each component of $\Sigma_0$, obtaining 
$$D^{N_i}_{F_i}:L_1^p\Gamma(N_i\otimes L(F_i))\to L^p\Omega^{0,1}(N_i\otimes L(F_i)).$$
The Chern class is diminished by $|F_i|$, thus we get
\medskip
\noindent{\bf Proposition 2.} \hskip3mm
{\it Under the hypotheses of Proposition 1, let 
$F\subset\Sigma_0$ be a finite subset such that
$$|F_i|+|df_i^{-1}(0)|<\langle c_1(f_i^*TV),\Sigma_i\rangle\hskip3mm\forall i=1,\cdots,r.$$ 
Then $D^{N_0}_F$ is onto.}
\bigskip
Combining this with Theorem 1', we get Corollary 2.
\vskip10mm
\noindent{\bf 7. $J$-curves and symplectic surfaces in dimension $4$}
\medskip
In this section we assume that $(V,\omega)$ is a compact symplectic $4$-manifold, 
so that there is a non-empty and contractible subspace ${\cal J}_\omega(V)\subset{\cal J}(V)$
of $\omega$-positive  almost complex structures.
\smallskip
We study the isotopy problem for symplectic surfaces, in relation
with $J$-holomorphic curves. 
\smallskip
Then we restrict to the case $(V,\omega)=(\C\P^2,\omega_0)$
and we give explicit sufficient conditions on a $J$-curve (resp. a curve and 
a finite subset)
to satisfy the hypotheses of Corollary 1 (resp. Corollary 2). and finally we apply Corollary 2 
to the case of surfaces of degree $3$.

\bigskip

\noindent{\bf 7.1. The isotopy problem for symplectic surfaces}
\medskip
Let \ $S\subset V$ be a symplectic surface, ie a real embedded surface such that 
$\omega_{|TS}$ never vanishes. Such a surface is connected and canonically oriented, and 
has a homology class $A$.

\medskip
The relation with $J$-holomorphic curves is the following fundamental observation of Gromov:
$S$ is symplectic if and only 
if there exists an almost complex structure $J\in {\cal J}_\omega(V)$ such that 
$S=f_0(\Sigma_0)$ where $[\Sigma_0,f_0]$ is 
an embedded $J$-holomorphic curve. Note that the genus of $S$ is $$g=g_a(A)={A.A-c_1(TV),A\over2}+1.$$
\medskip
Let us define $$\overline{\mkern -4mu{\cal M}}(V,J,A)=\overline{\mkern -4mu{\cal M}}_{g(A)}(V,J,A),$$
which we know is a compact space. It contains as an open subset the connected smooth curves:
$${\mkern -4mu{\cal M}}(V,J,A)={\mkern -4mu{\cal M}}_{g(A)}(V,J,A).$$
If $J$ is integrable, \ $\overline{\mkern -4mu{\cal M}}(V,J,A)$  sits over 
the space ${\cal D}(V,A)$ of divisors in the class $A$ (which is an algebraic variety), 
the projection being a homeomorphism from  ${\cal M}_{g(A)}(V,J,A)$ 
to the open subset ${\cal D}^s(V,A)$ of connected smooth curves. 

\smallskip
The space ${\cal D}^s(V,A)$ always has a finite number of connected
components (which are the same as path-connected components), thus one may hope that ${\cal M}_{g(A)}(V,J,A)$ also has a finite number
of path-connected components
for any fixed $J\in{\cal J}_\omega(V)$. More daringly, one may hope that it is also true for
for the space ${\cal S}(V,A)$ of connected symplectic surfaces in the class $A$ (which are necessarily of genus $g(A)$). In some cases,
 as for of $\C\P^2$, one may
even hope that ${\cal S}(V,A)$ is connected:

\medskip
\noindent{\bf Question.}\hskip3mm Let $d$ be a positive integer. Are any two symplectic surfaces of degree $d$ in
$\C\P^2$ symplectically isotopic ? 

Equivalently, is every symplectic surface in $\C\P^2$ symplectically isotopic
to an algebraic curve (necessarily smooth and of the same degree) ?

\medskip
In general however, it is not true that  this is not true that  ${\cal S}(V,A)$ has a finite number
of path-connected components, as shown recently by R. Fintushel and R. Stern [FS]: if $V$ is simply-connected and
${\cal S}(V,A)$
 contains a torus $T$
with zero self-intersection which can degenerate to a rational curve with a cusp, then there
are infinitely many tori $T_n$ in ${\cal S}(V,2A)$ which are pairwise not differentiably isotopic. In fact, all pairs 
$(V,T_n)$ are differentiably different. Note that the class $2A$ satisfies $\langle c_1(TV),2A\rangle=0$
 so that the results of this paper do not apply. This phenomenon is related to the possibility of describing 
a model of birth of $J$-holomorphic curves.
\medskip
The question above has a positive answer for $d=$ $1$ or $2$ by [G]. Note that these are the cases where
 the genus
is zero. Also, for $d=1$, an embedded sphere is always topologically isotopic to a complex line, and it is an open
question whether it is differentiably isotopic. For $d>1$, there are many surfaces of degree 
$g(d)={(d-1)(d-2)\over2}$, which are already topologically knotted, which makes the above question more
interesting.
\medskip
In paragraph 7.3 we shall recall Gromov's proof for $d=1$ and $2$, and explain why it does not directlly extend to $d\ge3$.
\bigskip
\noindent
{\bf 7.2. Hypothetical proof of the isotopy property}
\medskip
In this paragraph we give a positive answer to the above question in $\C\P^2$, depending on 
plausible properties of $J$-curves.
\medskip
Let $S\subset\C\P^2$ be a symplectic surface of degree $d$. One can
find an almost
complex structure $J\in {\cal J}_{\omega_0}$ such that $S$ is $J$-holomorphic. 
\smallskip
Let $(J_t)$  be a path in ${\cal
J}_{\omega_0}$ from the
standard structure to $J$. 
For $A=d[L]$, denote
$$\overline{\mkern -4mu{\cal M}}(t,d)=
\overline{\mkern -4mu{\cal M}}(\C\P^2,J_t,d)\,\,,\,\,{\cal M}(t,d)=
{\cal M}(\C\P^2,J_t,d),$$ 
 and denote  by \ $\overline{\mkern -4mu{\cal M}}^{\rm sing}(t,d)$ the
subset of singular
curves, ie $$\overline{\mkern -4mu{\cal M}}^{\rm sing}(t,d)=\overline{\mkern -4mu{\cal M}}(t,d)\setminus
{\cal M}(t,d).$$
Define also the extended spaces
$$\overline{\hbox{\cal M}}(d)
=\displaystyle\bigcup_{0\le t\le1}{\cal M}(t,d)\times\{t\}=\hbox{\cal M}(d)
\cup\overline{\hbox{\cal M}}^{\rm sing}(d),$$ and let $\pi:\overline{\hbox{\cal M}}\to[0,1]$ 
denote the second projection.
\medskip
By Gromov's compactness theorem and the automatic regularity, we have the following properties.
\medskip
\noindent{\bf Properties}
\medskip
1) $\overline{\hbox{\cal M}}(d)$ is compact.
\medskip
2) $\hbox{\cal M}(d)$ is a topological manifold of real dimension $d(d+3)+1$, with boundary 
$\hbox{\cal M}(d,0)\cup\hbox{\cal M}(d,1)$, and the restriction of $\pi$ is a topological submersion.
\medskip
Note that $\hbox{\cal M}(d,0)$ is the connected  space of all smooth algebraic curves of degree $d$.
\medskip
\noindent{\bf Ideal situation}
\medskip
Assume that {\it $\overline{\hbox{\cal M}}^{\rm sing}(d)$ never locally disconnects 
$\overline{\hbox{\cal M}}(d)$}, ie each point of $\overline{\hbox{\cal M}}^{\rm sing}(d)$ admits arbitrarily
small neighbourhoods $U_i$ in $\overline{\hbox{\cal M}}^{\rm sing}(d)$ such that 
$U_i\cap\hbox{\cal M}(d)$ is connected (and thus path-connected). This the case if it has the structure of a subcomplex
of real codimension $2$, which one can hope would follow from a suitable stratification result. 

\smallskip Under this assumption, the properties 1) and 2) above easily imply the path-connectedness of 
$\hbox{\cal M}(d)$ and thus the symplectic isotopy of $S$ to an algebraic curve.
\medskip
An important generalization of this method arises if we consider curves containing a fixed finite subset $F\subset S$ 
in the image.
We define in a similar manner $\overline{\cal M}(t,d;F)$,  ${\cal M}(t,d;F)$,
$\overline{\cal M}^{sing}(t,d;F)$, $\overline{\hbox{\cal M}}(d;F)$,... 
\medskip
Then Gromov's compactness theorem and the ``automatic genericity'' imply the
\medskip
\noindent{\bf Properties}
\smallskip
1') $\overline{\hbox{\cal M}}(d;F)$ is compact.
\smallskip
2') If $|F|<3d-1$, $\hbox{\cal M}(d;F)$ is a topological manifold of real dimension $d(d+3)+1-|F|$, with boundary 
$\hbox{\cal M}(d,0;F)\cup\hbox{\cal M}(d,1;F)$, and the restriction of $\pi$ is a topological submersion.

\bigskip
\noindent{\bf Generalized ideal situation}

\medskip
Assume the following:
\medskip
\item{(*)}
{\it For a suitable choice of $F$ and the path $(J_t)$,
$\overline{\hbox{\cal M}}^{\rm sing}(d;F)$ never locally disconnects 
$\overline{\hbox{\cal M}}(d)$.}
\medskip
Then again the properties 1') and 2') imply the path-connectedness of 
$\hbox{\cal M}(d;F)$ and thus the symplectic isotopy of $S$ to an algebraic curve.

\medskip
\bigskip
\noindent{\bf 7.3. Curves of degree at most $3$}
\medskip
When $d\le3$, the above hypothetical proof actually works. For $d=1$ or $2$ this is due to Gromov.
\medskip
\noindent{\bf Degree 1}
\medskip The key fact is that here there are no singular curves! Thus properties 1) and 2) imply that 
$\overline{\hbox{\cal M}}(1)$ is a compact $4$-dimensional manifold(with a natural smooth structure)  with boundary
$\overline{\hbox{\cal M}}(1,0)\cup\overline{\hbox{\cal M}}(1,1)$, and that the projection $\pi$ is a submersion 
of index $0$. Thus $\overline{\hbox{\cal M}}(1)$ is diffeomorphic to a product 
$\overline{\hbox{\cal M}}(1,0)\times[0,1]\approx\C\P^{2*}\times[0,1]$, 
which gives the desired symplectic isotopy of
the symplectic surface to a complex line.
\smallskip
\noindent{\it Remark. }\hskip3mm When one fixes the maximal number $3d-1=2$ points, one thus obtain the famous result
of Gromov that through any $2$ distinct points there passes a unique $J$-line.
\medskip
\noindent{\bf Degree $2$.}
\medskip Here there are singular curves, which are of two types:  a union of two distinct $J$-lines, and
a double $J$-line. The formal complex dimension of the first is $4$, and of the second
is $2$, since actually it is equivalent to the space of simple lines. 
Note however that as a subspace of
$\overline{\cal M}(2,J)$ it is of complex dimension $4$, being the blow up of the space of double lines in 
$\C\P^5$.
\smallskip
Thus when we fix 
the maximal number $3d-1=5$ points, the maximal complex formal dimension is $-1$, 
thus the real dimension is $-2$.
By a classical argument using Sard-Smale (see [MS] for details, and also [B]), this implies that for
a generic choice of $F$ and $(J_t)$, $\overline{\cal M}^{sing}(d,t)$ is always empty. Thus we conclude exactly as 
for degree $1$.

 \medskip
\noindent{\bf Degree $3$: proof of Theorem 3}
\medskip
Note that now the genus is $1$. Here we cannot fix enough
to avoid singular curves: 
indeed the ``stratum''  corresponding to rational curves with one double point, has 
formal complex dimension $8$, thus we would need to fix $9=3d$ points, spoiling property 1'). Considering curves
of degree $3$  passing through $9$ points is similar to considering
curves in $\C\P^2$ blown up at $9$ points: then $c_1(A)$ and $c_1(N)$ vanishes, and one can give a model 
of birth of $J$-holomorphic curves. What is more, the examples of Fintushel and Stern show that this can actually
lead to an infinite number of isotopy classes. 
\smallskip
Thus we fix only $8$ points. Using the results of [B], one can prove that all other strata than the one 
corresponding to rational curves with one double point have formal complex codimension at most $7$, and thus
will not occur for a generic choice of $F$ and $(J_t)$. The crucial case is that of rational
curves with one cusp, which has dimension $7$ as a space of curves parameterized by $\C\P^1$ (although, as a 
subspace of the space of stable curves of genus $1$, it consists of two strata of dimension $8$).

\medskip
Thus if $F$ and the path $(J_t)$ are chosen generically, , a curve in ${\cal M}^{\rm
sing}(3,t;F)$ is rational with a double point, and $F$ is contained in its smooth part. 
\medskip
Then Corollary 2 applies, thus:
\smallskip
i) $\overline{\hbox{\cal M}}(3;F)$ a topological manifold of dimension $3$
is a topological manifold, with boundary $\overline{\hbox{\cal M}}(3,0;F)\cup\overline{\hbox{\cal M}}(3,1;F)$ 
and the projection to $[0,1]$ is a submersion. 
\smallskip
ii)  $\overline{\hbox{\cal M}}^{\rm sing}(3;F)$ is a $1$-dimensional locally flat submersion.
\medskip
Thus the non-disconnecting property (*) is again satisfied, which implies Theorem 3.

\bigskip
\noindent{\bf Comments}\hskip3mm We have not been able to extend this method above degree $3$: already in degree
$4$, one cannot avoid curves with a cusp, which as stable curves include a constant component which usually 
forbids the normal genericity.

\vskip1cm
\centerline{\bf References}

\bigskip
\noindent{[B]} J.-F. Barraud, {\it Nodal symplectic spheres in $\C\P^2$ with positive
self-intersection}, 
Intern. Math. Res. Not. {\bf9} (1999), 495-508.
\medskip

\noindent{[F]} A. Floer, {\it Morse theory for Lagrangian intersections},
J. Diff.Geom.
{\bf28} (1988), 513-547.
\medskip
\noindent{[FO]} K. Fukaya and K. Ono, {\it Arnold conjecture and Gromov-Witten
invariant for
general symplectic manifolds}, Topology {\bf38} (1999), 933--1048.
\medskip
\noindent[FS] Ronald Fintushel, Ronald J. Stern, {\it Symplectic surfaces 
in a fixed homology class}, preprint SG/9902028, to appear in the Journal of Differential Geometry.
\medskip
\noindent{[G]} M. Gromov, {\it Pseudo holomorphic curves in symplectic
manifolds},
Invent. Math. {\bf82} (1985), 307-347.
\medskip
\noindent{[HLS]} H. Hofer, V. Lizan, J.-C. Sikorav, {\it On genericity for
holomorphic
curves in $4$-dimensional almost-complex manifolds}, J. Geom. Anal. {\bf7} (1998), 149-159.
\medskip
\noindent{[HM]} Joe Harris and Ian Morrison, {\it Moduli of curves}, Grad.
Texts in
Math. {\bf187}, Springer, 1998.
\medskip
\noindent{[IS]}, S. Ivashkovich and V. Schevchishin, {\it Structure of the moduli space in the 
neighbourhood of a cusp-curve and meromorphic hulls}, Invent. Math {\bf136}
(1999), 571-602.
\medskip
\noindent[KM] M. Kontsevich and Y. Manin, {\it Gromov-Witten classes,
 quantum cohomology and enumerative geometry }, Commun. Math. Phys. {\bf164} (1994), 525-562.
\medskip
\noindent[LL] T.J. Li and A. Liu, {\it Symplectic structure on ruled surfaces and a generalized adjunction
formula}, Math. Res. Lett. {\bf2} (1995), 453--471.
\medskip
\noindent{[LT]} J. Li and G. Tian, {\it Virtual moduli cycles and Gromov-Witten
invariants of general symplectic manifolds}, in: {\it Topics in symplectic
$4$-manifolds
(Irvine, 1996)}, Int. Press, 1998, 47-83.
\medskip
\noindent{[MS]} D. McDuff and D. Salamon, {\it $J$-holomorphic curves and quantum
cohomology}, Univ. Lect. Series {\bf6}, Amer. Math. Soc., 1994.
\medskip
\noindent{[MW]} M. Micallef and B. White, {\it The structure of branched points in minimal
surfaces and in pseudoholomorphic curves}, Ann. of Math. {\bf 139} (1994), 35-85.
\medskip
\noindent[RT1] Y. Ruan and G. Tian, {\it A mathematical theory of quantum cohomology}, J. Diff. Geom. {\bf42} (1995),
259-367.
\medskip
\noindent[RT2] Y. Ruan and G. Tian, {\it Higher genus symplectic invariants and sigma models coupled with gravity},
Invent. Math. {\bf130} (1997), 455-516.
259-367.
\medskip
\noindent{[Sie1]} B. Siebert, {\it Gromov-Witten invariants for general
symplectic manifolds}, preprint dg-ga\break 9608005.
\medskip
\noindent{[Sie2]} B. Siebert, {\it Symplectic Gromov-Witten invariants}, in: {\it New 
trends in algebraic geometry}, 
London Math. Soc. Lect. Notes {\bf264} (1999), Cambridge Univ. Press, 375-424.
\medskip
\noindent{[Sik1]}, J.-C. Sikorav, {\it Local properties of $J$-holomorphic curves}, in: {\it Holomorphic curves in symplectic geometry},
Audin and Lafontaine ed., Progress in Math. {\bf 117}, Birkh\"auser, 165-189.
\medskip
\noindent{[Sik2]}, J.-C. Sikorav, {\it Singularities of $J$-holomorphic
curves}, Math. Z. {\bf226} (1997), 359-373.

\vskip10mm
\noindent
Jean-Claude Sikorav

\noindent ENS Lyon, UMPA (UMR CNRS 5669)

\noindent 46, all\'ee d'Italie

\noindent
F-69364 Lyon cedex 07, FRANCE

\noindent
sikorav@umpa.ens-lyon.fr

\end